\documentclass[12pt,a4paper]{article}
\usepackage{graphicx,tikz-cd} % Required for inserting images
\usepackage{amssymb}
\usepackage{amsmath,amsthm}
\usepackage{latexsym, amsfonts,hyperref}
\allowdisplaybreaks
\newcommand{\be}{\begin{equation}}
\newcommand{\ee}{\end{equation}}
\newcommand{\bea}{\begin{eqnarray}}
\newcommand{\eea}{\end{eqnarray}}
\newcommand{\bean}{\begin{eqnarray*}}
\newcommand{\eean}{\end{eqnarray*}}
\newcommand{\brray}{\begin{array}}
\newcommand{\erray}{\end{array}}
\newcommand{\ben}{\begin{equation}{nonumber}}
\newcommand{\een}{\end{equation}{nonumber}}

\newtheorem{dfn}{Definition}[section]
\newtheorem{thm}[dfn]{Theorem}
\newtheorem{lmma}[dfn]{Lemma}
\newtheorem{ppsn}[dfn]{Proposition}
\newtheorem{crlre}[dfn]{Corollary}
\newtheorem{xmpl}[dfn]{Example}
\newtheorem{rmrk}[dfn]{Remark}
\newcommand{\bdfn}{\begin{dfn}}
\newcommand{\bthm}{\begin{thm}}
\newcommand{\blmma}{\begin{lmma}}
\newcommand{\bppsn}{\begin{ppsn}}
\newcommand{\bcrlre}{\begin{crlre}}
\newcommand{\bxmpl}{\begin{xmpl}}
\newcommand{\brmrk}{\begin{rmrk}}
\newcommand{\edfn}{\end{dfn}}
\newcommand{\ethm}{\end{thm}}
\newcommand{\elmma}{\end{lmma}}
\newcommand{\eppsn}{\end{ppsn}}
\newcommand{\ecrlre}{\end{crlre}}
\newcommand{\exmpl}{\end{xmpl}}
\newcommand{\ermrk}{\end{rmrk}}
\newcommand{\bbc}{\mathbb{C}}

\newcommand{\q}{\mathcal{Q}}

\newcommand{\clb}{{\cal B}}
\newcommand{\clc}{{\cal C}}

\newcommand{\clh}{{\cal H}}

\newcommand{\clk}{{\cal K}}

\def\a*{{\cal A}_{h,*}}
\def\B{{\cal B}(h)}
\def\B1{{\cal B}_1(h)}
\def\b{{\cal B}^{\rm s.a.}(h)}
\def\b1{{\cal B}^{\rm s.a.}_1(h)}

\newcommand{\ot}{\otimes}

\begin{document}
\[
\]
\begin{center}
{\large {\bf  Second invariant cohomology of some finite-dimensional Hopf algebras}}\\
by\\
{\large Debashish Goswami{\footnote{ Partially supported by JC Bose National Fellowship given by SERB, Govt. of India.}}
},
{\large Kiran Maity{\footnote{ Support by a research fellowship from Indian Statistical Institute is gratefully acknowledged.}}}\\
{ Stat-Math Unit, Kolkata,}\\
{ Indian Statistical Institute,}\\
{ 203, B. T. Road, Kolkata 700 108, India}\\
{e mails: debashish\_goswami@yahoo.co.in, maitykiran01@gmail.com}

\end{center}

\begin{abstract}
We use categorical description of the invariant 2-cohomology group of Hopf algebra to compute such cohomology for two finite dimensional Hopf algebras: the group ring of $Z_8\rtimes Aut(Z_8)$ and Kac-Paljutkin algebra. For the first of these two examples, our categorical approach helps to settle the problem of computing this cohomology, which was left open in \cite{dtwist}, where only some partial information about this cohomology was obtained.
\end{abstract}

 \section{Introduction}
Unitary dual 2-cocycles of a compact  group was first introduced by Landstand\cite{ldstd} and Wassermann \cite{wasse} to study full multiplicity ergodic actions. In the context of Hopf algera, it was used by Drinfeld in his work \cite{Drin} on quasi-Hopf algebra, which is known as Drinfeld twist. The idea of Drinfeld is  used by Etingof, Gelaki \cite{et1},\cite{et2} to construct and categorize the triangular semisimple Hopf algebra and to study the fiber functor on their representation categories. The classical theory of full multiplicity ergodic action of Wasserman was extended to compact quantum group setting by  Vaes, Bichon, An de Rijdt \cite{etn}. They proved that every dimension preserving unitary fiber functor on corepresentation category of a compact quantum group (CQG) corresponds to a unitary 2-cocycle of the dual of that compact quantum group.

In this paper, we mainly focus on invariant 2-cocycles on certain finite dimensional Hopf algebra. For a corepresentation category $Corep(\q)$ of compact quantum group $(\q,\Delta)$, the group of tensor structures (up to isomorphism) of identity functor  of $Corep(\q)$ is called the  second invariant cohomology group of the dual of $(\q,\Delta).$ It is the group of isomorphism classes of tensor autoequivalences of $Corep(\q)$ which act trivially on $Corep(\q).$ Kassel, Guillot \cite{dtwist} studied  invariant 2-cohomology for group algebra of finite groups. But their approach is not categorical . We use completely different approach based on categorical description to  compute invariant cohomology $H^{2}_{uinv}(\cdot, S^{1})$ and $H^{2}_{inv}(\cdot ,\bbc-\{0\})$
% we will explicitly compute 
 for two interesting finite dimensional Hopf algebras, namely  the group ring of $G:=Z_{8}\rtimes Aut(Z_{8})$ and the dual of Kac-Paljutkin algebra.  Let us mention that in \cite{dtwist}, partial results about $H^{2}_{inv}(C^{\ast}(G), S^{1})$ were obtained by quite tedious calculation using other methods and also using computer. In that paper, it could only be proved that the order of the cohomology group is either 2 or 4, without a definite conclusion. On the other hand, our categorical approach seems to simplify the computations to some extent and we are able to conclusively identify the group as $Z_2,$ without computer aided calculations. %alternative description of cohomology in term of unitary/invertibe fiber functors described in \cite{snbook}. First we will introduce  some theory of invariant 2-cocycle of DQG, which are mostly follows from \cite{classification of non kac type quantum group of SU(n),symmetric invariant cocycle}.

  The paper is organized as follows:  In section 2, we recall the definition of fiber functor on corepresentation category of a compact quantum group. Then we define the second invariant cohomology $H^{2}_{uinv}(\cdot, S^{1})$ and $H^{2}_{inv}(\cdot ,\bbc-\{0\})$ for a compact quantum group $\q.$ We proved some useful facts about invariant cohomology.

  In section 3 , we show that 
  \begin{align*}
      H^{2}_{uinv}(C^{*}(G), S^{1}) \cong H^{2}_{inv}(C^{*}(G) ,\bbc-\{0\}) \cong Z_2,
  \end{align*}
  %H^{2}_{uinv}(C^{*}(G), S^{1})$ and $H^{2}_{inv}(C^{*}(G) ,\bbc-\{0\})$ are isomorphic to $Z_2,$ 
    where $C^{*}(G)$ is the group algebra of  $Z_8\rtimes Aut(Z_8)$.

  In section 4, we prove that $H^{2}_{uinv}(\q_{kp}, S^{1})$ and $H^{2}_{inv}(\q_{kp} ,\bbc-\{0\})$ are trivial for the 8- dimensional non-cocommutative, non commutative Kac-Paljutkin algebra.
\section{Preliminaries}
Here we will not give the detailed definition of a compact quantum group (CQG) and rigid $C^{\ast}$ tensor category $\clc$. We will refer the reader to \cite{snbook} for it. For our purpose, we only require the following
\bxmpl 
The category of finite dimensional Hilbert space, denoted by $Hilb_{f}$, which has finite dimensional Hilbert spaces as the set of objects, morphisms given by linear maps from $\clh$ to $\clk.$ $\clh\otimes \clk$ is the usual tensor product between two Hilbert spaces. Associativity morphism $\alpha_{\clh_{1},\clh_{2},\clh_{3}}$ is the identity map from $(\clh_{1}\ot \clh_{2})\ot \clh_{3}$ to $\clh_{1}\ot (\clh_{2}\ot \clh_{3}).$ $1_{\bbc}$ is unit object of this category $Hilb_{f}.$ $Hilb_{f}$ is a strict rigid $C^{\ast}$ tensor category.
\exmpl
\bxmpl Representation category of a compact group is denoted by $Rep(G).$ Objects of this category are denoted by $(\pi,\clh_{\pi})$, where $\pi:G\to U(\clh_{\pi})$ is a finite dimensional unitary representation of $G$ on $\clh_{\pi}.$ Morphism between two objects $(\pi_{1},\clh_{\pi_{1}})$ $(\pi_{2},\clh_{\pi_{2}})$ is defined by 
   \begin{align*}
    Mor((\pi_{1},\clh_{\pi_{1}}),(\pi_{2},\clh_{\pi_{2}}))=\{T\in\clb(\clh_{\pi_{1}}, \clh_{\pi_{2}}):T\pi_{1}(g)=\pi_{2}T(g),~for~all~g\in G\}.
   \end{align*} Tensor product between two objects $(\pi_{1},\clh_{\pi_{1}})$ $(\pi_{2},\clh_{\pi_{2}})$ is defined by $(\pi_{1}\ot \pi_{2},\clh_{\pi_{1}}\ot \clh_{\pi_{2}}),$  where $\pi_{1}\ot \pi_{2}$ is the usual tensor product of two group representations and $\clh_{\pi_{1}}\otimes \clh_{\pi_{2}}$ is the usual tensor product between two Hilbert spaces. Associativity morphism is the identity map and unit object is the trivial representation on $\bbc.$ $Rep(G)$ is a strict $C^{\ast}$ tensor category.
\exmpl

\bxmpl
Unitary corepresentation category of a CQG denoted by $Corep(\q)$. Objects are given by $(U,\clh_{U}),$ where $U\in B(\clh_{U})\ot\q $ is a finite dimensional unitary corepresentation of $\q$ on $\clh_{U}$. Morphism between two objects $(U,\clh_{U}),(V,\clh_{V})$ is given by \begin{align*}
    Mor((U,\clh_{U}),(V,\clh_{V}))=\{T\in B(\clh_{U},\clh_{V}):(T\ot 1)U=V(T\ot 1)\}
\end{align*} Tensor product of two objects $(U,\clh_{U}),(V,\clh_{V})$ defined by $(U_{13}V_{23},\clh_{U}\ot \clh_{V}).$ Associativity morphism is the identity map and unit object is $Id_{\bbc}\ot 1_{\q}.$ For a projection $P\in End(U,\clh_{U})$, $((P\ot 1)U,P\clh_{U})$ is a subobject of $(U,\clh_{U}).$ $Corep(\q)$ is a strict $C^{\ast}$ tensor category.
\exmpl

\bdfn 
  A fiber functor $F:\clc\to Hilb_{f}$ is a tensor functor which   is   faithful and exact.
\edfn
\bxmpl
Let $F^{Nat}:Corep(\q)\to Hilb_{f}$ is a fiber functor defined by $F^{Nat}(U,\clh_{U})=\clh_{U}$, where $(U,\clh_{U})$ is an object of this category , identity on the morphisms and also $F_{2}(U,V),F_{0}$ are both identity maps.
\exmpl
 Now, we briefly discuss about second cohomology of Hopf *-algebra.  In this section, let us denote by $\ot$ the algebraic tensor product for Hopf * -algebra being considered.
 Let $(\q,\Delta)$ be a Hopf algebra and \begin{align*}
      &\Delta_{i}:\q^{\otimes n}\rightarrow \q^{\otimes (n+1)}~ such~that~\\
      &\Delta_{i}=id\otimes...\otimes\Delta\otimes...\otimes id,
      \end{align*} $\Delta$ is in the $i$-th position for $i\in\{1,2,.,n\}$ and we define $\Delta_{0}(x)=1\otimes x~and~\Delta_{n+1}(x)=x \otimes 1.$ So that $\Delta_{i}$ defined for $i=0,1,..,n+1.$ A left $n-$ cochain $\chi$ is an invertible element of $H^{\otimes n}.$  Coboundary  of a left $n$-cochain $\chi$ is a $(n+1)$-cochain\begin{align*}
          \delta_{\chi}=(\Pi_{i=0}^{even}\Delta_{i}(\chi))(\Pi_{i=1}^{odd}\Delta_{i}(\chi^{-1}))
      \end{align*}
    \bdfn 
    A left n-cochain $\chi\in \q^{\otimes n}$ is said to be a n-cocycle if $\delta_{\chi}=1.$ A left n-cocycle is said to be counital if $\epsilon_{i}\chi=1$, where $\epsilon_{i}=id\otimes...\epsilon\otimes...id$, $\epsilon$ at i-th position. 
    \edfn
    \bxmpl A left 1-cocycle $\chi$ is an invertible element in $\q$ such that $\Delta(\chi)=\chi\otimes\chi$ and it is automatically counital and any left 2-cocycle $\chi\in \q^{\otimes 2}$ is satisfy the equation\begin{align*}
        (1\otimes\chi)(id\otimes \Delta)(\chi)=(\chi\otimes 1)(\Delta\otimes id)(\chi)
    \end{align*} and it is counital if $(\epsilon\otimes id)(\chi)=1.$
    \exmpl
    \blmma Let $G$ be a finite group then any left 1-cocycle $\chi$ for the ring of continuous function $C(G)$ is a group homomorphism $\chi:G\to \bbc-\{0\}.$ A counital left 2-cocycle $\chi$ of $C(G)$ is a normalized  complex valued 2-cocycle of $G.$
    \elmma
    \vspace{2mm}
   \bdfn Let $\Omega_{1},\Omega_{2}$ be left 2-cocycles of $\q.$       $\Omega_{1},\Omega_{2}$ are said to be cohomologous if there exists an invertible element $h$ in $\q$ such that $\Omega_{2}=(h\otimes h)\Omega_{1}\Delta(h^{-1}).$
   \edfn
  Let $ {\bf H^{2}(\q,\bbc^{\ast})}$ be the set of cohomology classes of 2-cocycles. It does not necessarily form a group.
  \bdfn  A left 2-cocycle $\Omega$ is said to be unitary 2-cocycle if and only if $\Omega$ is a unitary element of $\q\otimes \q$ and two unitary left 2-cocycles are said to be unitarily cohomologous if and only if there exists a unitary element $u$ such that $\Omega_{2}=(u\otimes u)\Omega_{1}\Delta(u^{-1}).$
  \edfn
   Let $ {\bf H^{2}(\q,S^{1})}$ be the set of unitary cohomology classes of unitary left 2-cocycles.

    %Similarly one can define right cochain and right cocycles. 
     We can similarly define right $n-$cochain and $n-$cocycle. As we  only need right 1 and 2-  cocycle, let us define them below.
     \bdfn A right (unitary) 1-cocycle $\chi'$ is an invertible (a unitary) element in $\q$ such that $\Delta(\chi')=\chi'\ot\chi'$ and an invertible (a unitary) element $\chi'\in \q\ot \q,$ is said to be a (unitary) right 2-cocycle if it satisfying the equation 
     \begin{align*}
         (id\ot \Delta)(\chi')(1\ot \chi')=(\Delta\ot id)(\chi')(\chi'\ot 1).
     \end{align*}
     \edfn

     \brmrk $\Omega$ is  left 2-cocycle if and only if $\Omega^{\ast}$ is right 2-cocycle.
     \ermrk
    \bdfn 
    An element of $\q^{\otimes n}$ is said to be invariant if it commutes with the elements in the image of $\Delta^{n-1}:\q\rightarrow \q^{\otimes n},$ where $\Delta^{n-1}$ is defined inductively as follows:$\Delta^{1}=\Delta,$ and $\Delta^{k}$ is obtained by applying $\Delta$ to any of the factors of $\Delta^{k-1}.$
    \edfn
    \brmrk
     An invariant element of $\q\ot\q$ is a left cochain/cocycle if and only if it is also a right cochain/cocycle. Hence we will simply call them invariant cochain/cocycle.
     \ermrk
    \blmma 
    If $\Omega_{1},\Omega_{2}\in \q^{\otimes 2}$ are invariant unitary 2-cocycles, then $\Omega_{1}\Omega_{2}$ is a unitary 2-cocycle and $\Omega_{1}^{\ast}$,$\Omega_{2}^{\ast}$  are both unitary 2-cocycles.
    \elmma
     \vspace{2mm}
     Let $A^{1}(\q)$ be a set of all central invertible elements of $\q$ and $A^{2}(\q)$ be a set of all invariant invertible 2-cocycles of $\q.$
     \blmma $\delta:A^{1}(\q)\rightarrow A^{2}(\q)$ is a group homomorphism and image $\delta$ is a central subgroup of $A^{2}(\q).$
     \elmma
     \bdfn Invariant 1-cohomology and 2-cohomology of $\q$ is given by \begin{align*}
         &H_{inv}^{1}(\q,\bbc-\{0\})=ker(\delta:A^{1}(\q)\to A^{2}(\q)),\\
         &H_{inv}^{2}(\q,\bbc-\{0\})=A^{2}(\q)/Image(\delta).\\
         \end{align*}
         \edfn

         \bxmpl [Theorem 7.1 of \cite{dtwist}]\label{coexam}
         
 $H_{inv}^{2}(\q,\bbc-\{0\}) = 1$ if $\q=C^{*}(G)$ and $G$ belongs to the following list of
finite groups:

i) the simple groups,

(ii) the symmetric groups~$S_n$,

(iii) the groups $SL_n(F_q)$,

(iv) the groups $GL_n(F_q)$ when $n$ is coprime to $q-1$.
\exmpl

%\bxmpl[proposition 7.7 of \cite{dtwist}] $H_{inv}^{2}(\q,\bbc-\{0\}) =Z_{2}$, if $\q=\bbc^{\ast}(A_{4})$, where $A_{4}$ is the alternating group.
%\exmpl
%\blmma[Theorem 7.4 of \cite{etin}] If $G$ is a connected affine algebraic group then $H_{inv}^{2}(\bbc^{\ast}(G),\bbc^{\ast}) $ is a commutative group.
%\elmma
\brmrk 
For a finite group $G$, $H_{inv}^{2}(C^{\ast}(G),\bbc-\{0\})$ can be a noncommutative group.
\ermrk

Let $A_{u}^{1}(\q)$ be a set of all central unitary elements of $\q$ and $A_{u}^{2}(\q)$ be a set of all invariant unitary 2-cocycles of $\q.$ Similarly, we can prove that $\delta:A^{1}_{u}\rightarrow A^{2}_{u}$ is a group homomorphism and $\delta(A^{1}_{u})$ is a central subgroup of $A^{2}_{u}(\q).$
\bdfn Unitary Invariant 1-cohomology and 2-cohomology of $\q$ are given by \begin{align*}
         &H_{uinv}^{1}(\q,s^{1})=ker(\delta:A_{u}^{1}(\q)\to A_{u}^{2}(\q),\\
         &H_{uinv}^{2}(\q,s^{1})=A_{u}^{2}(\q)/\delta(A_{u}^{1}(\q)).\\
         \end{align*}
         \edfn
 The following lemma follows from lemma (3.1.5) of \cite{snbook}.
\blmma There exists an injective group homomorphism $\theta$ , is given by \begin{align*}
    &\theta:H_{uinv}^{2}(\q,s^{1})\rightarrow H_{inv}^{2}(\q,\bbc-\{0\})\\
    &\theta([\Omega])=[\Omega],
\end{align*} where $\Omega$ is an invariant unitary 2-cocycle in $\q.$
\elmma
 %We can easily extend the notion of invertible and unitary left/right/invariant cocycles to both CQG and DQG and also von Neumann bialgebra, by replacing the algebraic tensor product by suitable $C^{\ast}$ or von Neumann algebraic tensor product, and also replaces $\q\ot\q$ by $M(\q\ot \q)$ for a DQG $\q.$ Counital cocycles will make sense under some conditions, we have the counit is bounded e.g for DQG.

Let  $(\q,\Delta)$ be a finite dimensional CQG, $I$ be the set of mutually inequivalent irreducible corepresentations. It follows from \cite{etn} that there is a fiber functor  $\phi_{\Omega}$ on $Corep(\q)$ such that \begin{align}
    H_{\phi_{\Omega}(x)}=H_{x},~~\phi_{\Omega}(S)=\Omega^{-1}S~~\phi_{\Omega}(T)=\Omega^{-1}_{2}T,
\end{align} where $\Omega$ is an invertible 2-cocycle of dual Hopf *-algebra of $\q,$  $,x,y,z\in I$ and for all $S\in Mor(x,y\ot z),T\in Mor(a,x\ot y\ot  z).$ $\phi_{\Omega}$ defines a new dual Hopf *-algebra $\widehat{\q_{\Omega}}$, where the * algebra $End(\phi_{\Omega})=End(F_{Nat})=\Pi_{x\in I} B(\clh_x)$ and coproduct is given by $\hat{\Delta}_{\Omega}(a)=\Omega\hat{\Delta}\Omega^{\ast}.$ $\phi_{\Omega}$ is a unitary monoidal equivalence between $Corep(\q)$ and $Corep(\q_{\Omega}).$

If $\Omega$ is an invariant 2-cocycle then $\hat{\Delta}_{\Omega}(a)=\hat{\Delta}(a)$ for $a\in \Pi_{x\in I} B(\clh_{x}).$   Algebra structure of $\hat{\q}_{\Omega}$ same as $\hat{\q}.$ For further details on the algebraic properties of dual deformed CQG, we refer the reader to \cite{BG},\cite{Gj}.

\blmma $\phi_{\Omega}$ is an a monoidal autoquivalence of $Corep(\q)$, where $\Omega$ is an invertible invariant 2-cocycle. In case $\Omega$ is unitary, this gives unitary monoidal autoequivalence.
\elmma
\begin{proof}
   $\Omega^{-1}$ is also a unitary invariant 2-cocycle of dual Hopf *- algebra of $\q.$ By choosing the fiber functor $\phi_{\Omega^{-1}},$ it follows that $\phi_{\Omega}\phi_{\Omega^{-1}}=\phi_{\Omega^{-1}}\phi_{\Omega}=Id_{corep(\q)}.$ 
\end{proof}

Our approach will be based on the following result, contained in chapter (3.1) of \cite{snbook}.
%$H^{2}_{uinv}(\widehat{\q},\hat{\Delta})$ is the group of all monoidal autoequivalences of $Corep(\q)$ that are naturally isomorphic to identity functor( tensor structure may not br trivial), but considered up to monoidal isomorphism(Chapter (3.1) of \cite{snbook}).
For the sake of completeness we give an outline of the proof here as well.
\bthm $H^{2}_{uinv}(\widehat{\q},\hat{\Delta})$ is isomorphic with the group of unitary isomorphism class of all unitary monoidal autoequivalences of $Corep(\q)\cong Rep(\hat{\q})$ which are naturally isomorphic to the identity functor. A similar statement holds  for $H^{2}_{inv}(\hat{\q}, \bbc-\{0\})$ without the requirement of unitary.
\ethm
\begin{proof}
Let us prove the case of unitary cohomology only, as the other case is very similar. Given a unitary invariant 2-cocycle $\Omega,$ the UTF $\phi_{\Omega}$ clearly gives a unitary monoidal autoequivalence of $Corep(\q),$ identifying $Corep(\q_{\Omega})$ with $Corep(\q).$

To prove the converse, assume that $\phi$ is an autoequivalence with the stated properties. As $\phi(x)\cong x$ for every object $x$ of $Corep(\q),$ $\phi$ is dimension preserving, hence it must be isomorphic with one of the form $\phi_{\Omega}$ for some unitary right 2-cocycle. As $\phi_{\Omega}$ is isomorphic with the identity functor, for each irreducible object $x,y,z$, we have unitary morphisms $\eta_{x},\eta_{y\ot z}$
such that $\phi_{\Omega}(t)=\eta^{-1}_{y\ot z} t\eta_x$ for all $t\in Mor(x, y\ot z).$ Moreover, as $x$ is irreducible, $\eta_x$ is a nonzero constant multiple of identity, say $c_xI_{x}.$ This implies that the component $\Omega_{y\ot z}\in Mor(y\ot z, y\ot z)$ is given by $\eta_{y\otimes z}\cdot c^{-1}_x,$ so it commutes with $\hat{\Delta}(\Theta)_{y,z}=\Theta_{y\ot z}$ for any $\Theta$ in $End(F_{Nat}).$ This proves that $\Omega$ is an invariant 2-cocycle.

\end{proof}
\brmrk Note that any monoidal functor $\phi: Corep(\q)\to Corep(\q)$ is natuarally isomorphic with the identity functor if and only if $\phi(x)\cong x$ for all objects $x$ of $Corep(\q).$
\ermrk
%Let $\psi$ be a UTF that isomorphic to identity functor.
%For each irreducible objects $x,y$ of $corep(\q)$ there exists a natural map $\eta_{x},\eta_{y}$ such that $\eta_{x}$ is an isomorphism between $\psi(x)$ and $x, $ $\psi(x)= (\eta^{-1}_x\ot 1) x (\eta_{x}\ot 1).$ Hence $\psi$ is a dimension preserving fiber functor. Then there exists a 2-cocycle $\Omega$ such that $\psi\equiv \phi_{\Omega}.$ Without loss of generality, we assume that $\psi=\phi_{\Omega}.$

% It is already known that $\phi_{\Omega}$ is a monoidal equivalence between $Corep(\q)$ and $Corep(\q_{\Omega}).$ It is also given that $\phi_{\Omega}$ is an autoequivalnce of $Corep(\q).$ 
 %If $E'$ is an equivalence between a UTC $C$ and $Corep(\q')$ and $F'$ is the canonical fiber functor on $Corep(\q').$ If $F$ is a fiber functor on $C$ and $F\cong F'E'$ then 
 %\begin{align}
 %    (End(F),\delta_{F})\cong (End(F'),\delta_{F'})
% \end{align}

%If we chose $C=Corep(\q), (\q',\Delta)=(\q,\Delta),$ $E'=\phi_{\Omega}$ and $F=\phi_{\Omega}$ then $(\q_{\Omega},\Delta)=(\q, \Delta).$ So $\Omega$ is an invariant 2-cocycle of $\hat{\q_{\Omega}}=\hat{\q}.$
%$End(F)\cong \Pi_{x\in I} \eta_{x}$ and natural isomorphism given by $\Pi_{x\in I}.$

\section{Cohomology of the group ring {$Z_8\rtimes Aut(Z_8)$}{}} 
The group $G:=Z_8\rtimes Aut(Z_8)$ was  considered by G. E. Wall in his paper \cite{holg}, $Z_8\rtimes Aut(Z_8)$ is generated by $s,t,u$, where $s,t,u$ satisfies the relations 
\begin{align}
s^2=t^2=u^8=1, ~~~st=ts,~~sus^{-1}=u^3,~~~tut^{-1}=u^5~.
\end{align}
Let $\chi_{ijk}$ be characters of $G$, defined by \begin{align}
    \chi_{ijk}(u)=(-1)^k,~~\chi_{ijk}(s)=(-1)^i,~~\chi_{ijk}(t)=(-1)^j~~~i,j,k\in Z_2.
\end{align}

Let $\pi_{2}$ and $\pi^{'}_{2}$ be irreducible representations of $G$  on $\clh_{\pi_{2}}$ and $\clh_{\pi^{'}_{2}}$, where $\{e_1,e_2\}$ is an orthonormal basis for $\clh_{\pi_{2}}$ and $\{f_{1},f_{2}\}$ is an orthonormal basis for $\clh_{\pi^{'}_{2}}.$ $\pi_{2}$ is given by
\begin{align}
    \pi_{2}(u)=\begin{pmatrix}
\omega^2 & 0 \\
0 & \omega^6
\end{pmatrix},
\pi_{2}(s)=
\begin{pmatrix}
   0 & 1\\
   1 & 0
\end{pmatrix},
\pi_{2}(t)=Id_{\clh_{\pi_{2}}},
\end{align} where $\omega=e^{2\pi i/8}.$

$\pi^{'}_{2}$ is given by 
\begin{align}
  \pi^{'}_{2}(u)=\begin{pmatrix}
\omega^2 & 0 \\
0 & \omega^6
\end{pmatrix},
\pi^{'}_{2}(s)=
\begin{pmatrix}
   0 & 1\\
   1 & 0
\end{pmatrix},
\pi^{'}_{2}(t)= -Id_{\clh_{\pi^{'}_{2}}}.
\end{align}

Let $\pi_{4}$ be an irreducible representation of $Z_8\rtimes Aut(Z_8)$ on $\clh_{\pi_{4}}$ and assume that $\{x_{1},x_{2},x_{3},x_{4}\}$ is an orthonormal basis for $\clh_{\pi_{4}}.$ $\pi_{4}$ is given by
\begin{align}
    \pi_{4}(u)=\begin{pmatrix}
        \omega & 0 & 0 & 0\\
         0     & \omega^{3} & 0 & 0\\
         0     &  0 & \omega^5 & 0\\
         0     &    0        & 0  & \omega^7
    \end{pmatrix},
\pi_{4}(s)=\begin{pmatrix}
        0 & 1 & 0 & 0\\
         1     & 0 & 0 & 0\\
         0     &  0 & 0 & 1\\
         0     &    0        & 1  & 0
    \end{pmatrix},
    \pi_{4}(t)=\begin{pmatrix}
        0 & 0 & 1& 0\\
         0     & 0 & 0 & 1\\
         1     &  0 & 0 & 0\\
         0     &    1        & 0  & 0
    \end{pmatrix}.
\end{align}
By an easy calculation one can verify that all the above are mutually inequivalent irreducible representations of $G$ and looking at the dimension of $C^{\ast}(G),$ we deduce that they exhaust the all irreducible representations of $G.$

Let $I:=\{\chi_{ijk},\pi_{2},\pi^{'}_{2},\pi_{4}\}$ be the collection of all pairwise non-equivalent irreducible representations of $G.$
\blmma The following fusion rules hold:
\begin{enumerate}
    \item [1)]$\chi_{i_1j_1k_1}\ot \chi_{i_2j_2k_2}=\chi_{(i_1+i_2)(j_1+j_2)(k_1+k_2)}.$
    \item[2)] $\pi_{2}\ot \pi_{2}=\bigoplus \chi_{i0k}.$
    \item[3)]$\pi^{'}_{2}\ot \pi^{'}_{2}=\bigoplus \chi_{i0k}.$
    \item[4)]$\pi_{2}\ot \pi^{'}_{2}=\bigoplus \chi_{i1k}.$
\end{enumerate}\elmma
\begin{proof}
    Proof is omitted, as it is a straightforward verification.
\end{proof}

Let $\phi$ be a unitary autoequivalence  on $Corep(C(G))\cong Rep(G)$, where $C(G)$ is the CQG of all continuous functions of $G.$ Then $\phi$ is automatically a dimension-preserving fiber functor.

Let  $Mor(a,b\ot c\ot d):=\Upsilon^{b\ot c\ot d}_a$ define a basis of $Mor(a, b\ot c\ot d)$ of unit norm, where $a,b,c,d \in \{\chi_{ijk},\pi_2,\pi'_2\}.$ Assume that $\phi(Mor(a,b\ot c\ot d)):=\widetilde{\Upsilon}^{b\ot c\ot d}_a.$

Now, we will introduce some notation to simplify our mathematical expressions.\begin{align}
    &\widetilde{\Upsilon}^{\chi_{1}\ot \chi_{2}}_{\chi_{1}\chi_{2}}=\sigma(\chi_1,\chi_2)\Upsilon^{\chi_{1}\ot \chi_{2}}_{\chi_{1}\chi_{2}},  ~where~\chi_1,\chi_2~ are~ characters.\\
   &\widetilde{\Upsilon}^{\chi\ot\pi_{2}}_{\pi_2}=c_{\chi}\Upsilon^{\chi\ot\pi_2}_{\pi_2},~where~\chi\in\{\chi_{i0k}\}.\\ 
    &\widetilde{\Upsilon}^{\chi\ot\pi^{'}_{2}}_{\pi^{'}_2}=c^{'}_{\chi}\Upsilon^{\chi\ot\pi^{'}_2}_{\pi^{'}_2},~where~\chi\in\{\chi_{i0k}\}.\\
    &\widetilde{\Upsilon}^{\chi^{'}\ot\pi_{2}}_{\pi^{'}_2}=d_{\chi^{'}}\Upsilon^{\chi^{'}\ot\pi_2}_{\pi^{'}_2}, ~where~\chi^{'}\in \{\chi_{i1k}\}.\\
    &\widetilde{\Upsilon}^{\chi^{'}\ot\pi^{'}_{2}}_{\pi_2}=e_{\chi^{'}}\Upsilon^{\chi^{'}\ot\pi^{'}_2}_{\pi_2},~where~\chi^{'}\in\{\chi_{i1k}\}.\\
    &\widetilde{\Upsilon}^{\pi_{2}\ot\pi_{2}}_{\chi}=\lambda_{\chi}\Upsilon^{\pi_{2}\ot\pi_{2}}_{\chi},~where~\chi \in \{\chi_{i0k}\}.\\
    &\widetilde{\Upsilon}^{\pi^{'}_{2}\ot\pi^{'}_{2}}_{\chi}=\lambda^{'}_{\chi}\Upsilon^{\pi^{'}_{2}\ot\pi^{'}_{2}}_\chi,~where~\chi \in \{\chi_{i0k}\}.\\
    &\widetilde{\Upsilon}^{\pi_{2}\ot\pi^{'}_{2}}_{\chi^{'}}=\eta_{\chi^{'}}\Upsilon^{\pi_{2}\ot\pi^{'}_{2}}_{\chi^{'}},~where~\chi^{'}\in \{\chi_{i1k}\}.
   \end{align} 
    
%    Here we introduce the $6j$ symbol, where $a,b,c\in I,$

   % \begin{tikzpicture}
%\draw (-2,2) to (0,0) to (2,2);
%\draw (0,0)  to (0,-2);
%\draw node at (-2,2)     {$a$};
%\draw node at (0,-2)  {$c$};
%\draw node at (2,2) {$b$};
%\end{tikzpicture} $\in  Mor(c,a\ot b).$

%From the remark (1.4.20), we can conclude that the collection of all 

%\begin{tikzpicture}
%\%draw (-1,1) to (0,0) to (1,1);
%\d%raw (0,0)  to (1,-1);
%\dr%aw (2,0)  to  (2,1);
%\dra%w (2,0) to  (1,-1);
%\dra%w (2,-1) to  (1,-2);
%\%draw %node at (-1,1)     {$a$};
%\draw n%ode at (1,1)  {$b$};
%\draw node at (2,1) {$x$};
%\draw node at (1,-1) {$y$};
%\%draw node at (0,0) {$c$};
%\draw node at (2,0) {$x$};
%\end{tikzpicture} 
%is the same as the collection of all

%\begin{tikzpicture}
%\draw (-1,1) to (0,0) to (1,1);
%
%
%
%
%
%
%
%
%
%
%\draw (0,0)  to (-1,-1);
%\draw (-2,1)  to  (-2,0);
%\draw (-2,0) to  (-1,-1);
%\draw (2,-1) to  (1,-2);
%\draw node at (-2,1)     {$a$};
%\%draw node at (-1,1)  {$b$};
%\draw node at (1,1) {$x$};
%\draw node at (1,-1) {$y$};
%\draw node at (0,0) {$d$};
%\draw node at (2,0) {$x$};
%\end{tikzpicture}
Also assume that for $\chi=\chi_{i0k},$ $\chi'=\chi_{i1k}.$

By our defining notations, we can say that for any $T\in Mor(a, b\ot c)$, $\widetilde{T}=c_{0}T$ (for a suitable choice of $c_{0}\in\bbc$ that satisfies our predefined notations). From now on, we only write $\widetilde{\Upsilon}^{b\ot c}_{a}=\widetilde{T}=c_{0}$ instead of $c_{0}T$, where we fix a $T\in Mor(a,b\ot c),  ~and~a,b,c\in\{\chi_{ijk},\pi_{2},\pi^{'}_{2}\}.$

\blmma \label{le3.2} $\sigma$ is a 2-cycle on the group $\{\chi_{ijk}\}.$ 
\elmma
\begin{proof}
 Let $\chi_1,\chi_2,\chi_3$ be characters of $Z_8\rtimes Aut(Z_8).$
  From this diagram
\[\begin{tikzcd}
	{\chi_{1}\otimes \chi_{2} \otimes \chi_3} && {\chi_{1}\chi_{2}\otimes \chi_{3}} \\
	{\chi_{1}\otimes\chi_{2}\chi_{3}} && {\chi_{1}\chi_{2}\chi_{3}}
	\arrow["{\sigma^{\ast}(\chi_{1},\chi_{2})\otimes Id_{\chi_{3}}}", from=1-1, to=1-3]
	\arrow["{id_{\chi_{1}}\otimes\sigma^{\ast}(\chi_2,\chi_3)}"', from=1-1, to=2-1]
	\arrow["{\sigma^{\ast}(\chi_{1}\chi_{2},\chi_{3})}", from=1-3, to=2-3]
	\arrow["{\sigma^{\ast}(\chi_{1},\chi_{2}\chi_{3})}", from=2-1, to=2-3]
\end{tikzcd}\] we can say that
\begin{align*}
    \sigma(\chi_{1},\chi_{2}\chi_{3}) \sigma(\chi_{2},\chi_{3})=\sigma(\chi_{1}\chi_{2},\chi_{3})\sigma(\chi_{1},\chi_{2}).
\end{align*} So, $\sigma$ is a 2-cocycle on the group $\{\chi_{ijk}\}\cong Z_{2} \times Z_{2} \times Z_{2}.$ Without loss of generality, we assume that $\sigma$ is a normalized 2-cocycle.
\end{proof}
\blmma
$c^{2}_{\chi}=(c'_\chi)^{2}=1$ and $c_{\chi_{000}}=c'_{\chi_{000}}=1.$
\elmma
\begin{proof}
    For any  character $\chi,$ $\chi\ot\chi\equiv \chi_{000}$ and $Mor(\pi_{2},\chi\ot\chi \ot \pi_{2})\equiv\bbc .$ If we choose a $T\in Mor(\pi_{2},\chi\ot\pi_{2})$ then  $(Id_{\chi}\ot T)T\in Mor(\pi_{2},\chi\ot\chi \ot \pi_{2}).$ From equation (8),  we can say that $\widetilde{T}=c_{\chi}T.$ Now, we can observe that $(Id_{\chi}\ot \widetilde{T})\widetilde{T}=c^{2}_{\chi}(Id_{\chi}\ot T)T$. Hence $c_{\chi}^{2}=1$ as $\sigma(\chi,\chi)=1$ %$\phi(S)=S, ~for~S\in Mor(\pi_{2},\chi\ot\chi \ot \pi_{2})$
 and also it is easy to observe that $c_{\chi_{000}}=1.$

 Similarly, we can conclude that $(c^{'}_{\chi})^{2}=1$ and $(c^{'}_{\chi})^{2}=1.$
\end{proof}
\blmma
$c'_{\chi}=\psi(\chi)c_{\chi},$ where $\chi\in\{\chi_{i0k}\}$  and~ $\psi$ is a character on $\{\chi_{i0k}\}.$
\elmma
\begin{proof}
Let us assume that $\chi_{a},\chi_{b}\in\{\chi_{i0k}\}.$ Dimension of $Mor(\pi_{2},\pi_{2}\ot\chi_{a}\ot \chi_{b})=1$ .
 From this commutative diagram below
\[\begin{tikzcd}
	{\pi_{2}\otimes\chi_{a}\otimes\chi_{b}} && {} & {\pi_{2}\otimes\chi_{a} \chi_{b}} \\
	{\pi_{2}\otimes\chi_{b}} && {} & {\pi_{2}}
	\arrow["{ \sigma^{\ast}(\chi_{a},\chi_{b})}", from=1-1, to=1-4]
	\arrow["{c^{\ast}_{\chi_{a}}}"', from=1-1, to=2-1]
	\arrow["{c^{\ast}_{ab}}", from=1-4, to=2-4]
	\arrow["{c^{\ast}_{\chi_{b}}}"', from=2-1, to=2-4],
\end{tikzcd}\] we can conclude that
\begin{align}
    \label{15}\sigma(\chi_{a},\chi_{b})=\frac{c_{\chi_{a}}c_{\chi_{b}}}{c_{\chi_{a}\chi_{b}}}.
\end{align}

  From the following commutative diagram,
 \[\begin{tikzcd}
	{\chi_{a}\otimes\chi_{b}\otimes\pi^{'}_{2}} &&& {\chi_{a}\chi_{b}\otimes \pi^{'}_{2}} \\
	{\chi_{a}\otimes\pi^{'}_{2}} &&& {\pi^{'}_{2}}
	\arrow["{(\sigma(\chi_{a},\chi_{b}))^{\ast}}", from=1-1, to=1-4]
	\arrow["{(c^{'}_{\chi_{b}})^{\ast}}"', from=1-1, to=2-1]
	\arrow["{(c^{'}_{\chi_{a}\chi_{b}})^{\ast}}", from=1-4, to=2-4]
	\arrow["{(c^{'}_{\chi_{a}})^{\ast}}"', from=2-1, to=2-4],
\end{tikzcd}\] we will get 
\begin{align}
    \sigma(\chi_{a},\chi_{b})=\frac{c^{'}_{\chi_{a}}c^{'}_{\chi_{b}}}{c^{'}_{\chi_{a}\chi_{b}}}.
\end{align}
So, $\sigma(\chi_{a},\chi_{b})=\frac{c^{'}_{\chi_{a}}c^{'}_{\chi_{b}}}{c^{'}_{\chi_{a}\chi_{b}}}=\frac{c_{\chi_{a}}c_{\chi_{b}}}{c_{\chi_{a}\chi_{b}}}$ and this implies that $c^{'}_{\chi_{a}}=\psi'(\chi_{a})c_{\chi_{a}}$, where $\psi'$ is a character on the group $\{\chi_{iok}\}\cong Z_2\times Z_2.$  
\end{proof}
\blmma 
$e_{\chi'}d_{\chi'}=1$, where $\chi'\in\{\chi_{i1k}\}.$
\elmma
\begin{proof}
Let us assume that $\chi^{'}_{a},\chi^{'}_{b}\in \{\chi_{i1k}\}.$ Now, from the diagram below,
\[\begin{tikzcd}
	{\chi^{'}_{a}\otimes\chi^{'}_{b}\otimes\pi_{2}} &&& {\chi_{a} \chi_{b}\otimes \pi_{2}} \\
	{\chi^{'}_{a}\otimes\pi^{'}_{2}} &&& {\pi_{2}}
	\arrow["{(\sigma(\chi^{'}_{a},\chi^{'}_{b}))^{\ast}}", from=1-1, to=1-4]
	\arrow["{(d_{\chi^{'}_{b}})^{\ast}}"', from=1-1, to=2-1]
	\arrow["{(c_{\chi_{a}\chi_{b}})^{\ast}}", from=1-4, to=2-4]
	\arrow["{(e_{\chi^{'}_{a}})^{\ast}}"', from=2-1, to=2-4]
\end{tikzcd}\]
we will get \begin{align}\label{17}
\sigma(\chi^{'}_{a},\chi^{'}_{b})=\frac{d_{\chi^{'}_{b}}e_{\chi^{'}_{a}}}{c_{\chi_{a}\chi_{b}}}.   
\end{align}  Taking $\chi^{'}_{a}=\chi^{'}_{b}$ we have $d_{\chi^{'}_{a}}e_{\chi^{'}_{a}}=1.$
\end{proof}
\blmma \label{lemma3.6}

$e^{2}_{\chi'_{a}}=e^{2}_{\chi_{010}}\psi'(\chi_{a})$,where $\chi'_{a}\in\{\chi_{i1k}\}$ and assume  $\chi'_a=\chi_{i1k}$ when $\chi_a=\chi_{i0k}$.
\elmma
\begin{proof}
 We can obtain
\begin{align}\label{18}
\sigma(\chi^{'}_{a},\chi^{'}_{b})=\frac{e_{\chi^{'}_{b}}d_{\chi^{'}_{a}}}{c^{'}_{\chi_{a}\chi_{b}}}   
\end{align}
from the diagram below
\[\begin{tikzcd}
	{\chi^{'}_{a}\otimes\chi^{'}_{b}\otimes\pi^{'}_{2}} &&& {\chi_{a} \chi_{b}\otimes \pi^{'}_{2}} \\
	{\chi^{'}_{a}\otimes\pi_{2}} &&& {\pi^{'}_{2}}
	\arrow["{(\sigma(\chi^{'}_{a},\chi^{'}_{b}))^{\ast}}", from=1-1, to=1-4]
	\arrow["{(e_{\chi^{'}_{b}})^{\ast}}"', from=1-1, to=2-1]
	\arrow["{(c^{'}_{\chi_{a}\chi_{b}})^{\ast}}", from=1-4, to=2-4]
	\arrow["{(d_{\chi^{'}_{a}})^{\ast}}"', from=2-1, to=2-4].
\end{tikzcd}\]
 Comparing equations (\ref{17}) and (\ref{18}),  we get
\begin{align}\label{20}
    &\frac{e_{\chi^{'}_{b}}d_{\chi^{'}_{a}}}{c^{'}_{\chi_{a}\chi_{b}}}= \frac{e_{\chi^{'}_{b}}d_{\chi^{'}_{a}}}{\psi'(\chi_{a}\chi_{b})c_{\chi_{a}\chi_{b}}}=\sigma(\chi^{'}_{a},\chi^{'}_{b})=\frac{d_{\chi^{'}_{b}}e_{\chi^{'}_{a}}}{c_{\chi_{a}\chi_{b}}},\\
    &d_{\chi^{'}_{b}}e_{\chi^{'}_{a}}=\psi'(\chi_{a}\chi_{b})e_{\chi^{'}_{b}}d_{\chi^{'}_{a}}.
    \end{align}
    If $\chi^{'}_{b}=\chi_{010}$, we have $e^{2}_{\chi^{'}_{a}}=e^{2}_{\chi_{010}}\psi'(\chi_{a}).$
\end{proof}
\blmma $\lambda^{'}_{\chi_{b}}=\psi'(\chi_{b})c_{\chi_{b}}\lambda^{'}_{\chi_{000}}.$
\elmma
\begin{proof}
  
    \[\begin{tikzcd}
	{\chi_{b}\otimes\pi^{'}_{2}\otimes \pi^{'}_{2}} &&& {\pi^{'}_{2}\otimes\pi^{'}_{2}} \\
	{\chi_{b}\otimes \chi_{a}} &&& {\chi_{a}\chi_{b}}
	\arrow["{(c^{'}_{\chi_{b}}})^{\ast}", from=1-1, to=1-4]
	\arrow["{(\lambda^{'}_{\chi_{a}}})^{\ast}"', from=1-1, to=2-1]
	\arrow["{(\lambda^{'}_{\chi_{a}\chi_{b}}})^{\ast}", from=1-4, to=2-4]
	\arrow["{(\sigma(\chi_{a},\chi_{b})})^{\ast}"', from=2-1, to=2-4]
\end{tikzcd}\] From this diagram, we  observe that 
\begin{align}\label{21}
\sigma(\chi_{a},\chi_{b})=\frac{c^{'}_{\chi_{b}}\lambda^{'}_{\chi_{b}\chi_{a}}}{\lambda^{'}_{\chi_{a}}}=\frac{\psi'(\chi_{b})c_{\chi_{b}}\lambda^{'}_{\chi_{b}\chi_{a}}}{\lambda^{'}_{\chi_{a}}}.
\end{align}
After comparing equations (\ref{15}) and (\ref{21}), the following relation will occur
\begin{align}
   \lambda^{'}_{\chi_{ab}} c_{\chi_{a}\chi_{b}}=\psi'(\chi_{b})c_{\chi_{a}}\lambda^{'}_{\chi_{a}}.
\end{align}
Let us  assume that $\chi_{a}=\chi_{000}$. From equation (22), we will get \begin{align}
    \lambda^{'}_{\chi_{b}}=\psi'(\chi_{b})c_{\chi_{b}}\lambda^{'}_{\chi_{000}}.
\end{align}
  
\end{proof}

\vspace{4mm}

\blmma $\lambda_{\chi}=c_{\chi}\lambda_{\chi_{000}}.$
\elmma
\begin{proof}
    From the diagram below
\[\begin{tikzcd}
	{\chi_{a}\otimes\pi_{2}\otimes\pi_{2}} &&& {\chi_{a}\otimes \chi_{b}} \\
	{\pi_{2}\otimes\pi_{2}} &&& {\chi_{a}\chi_{b}}
	\arrow["{(\lambda_{\chi_{b}})^{\ast}}", from=1-1, to=1-4]
	\arrow["{(c_{\chi_{a}})^{\ast}}"', from=1-1, to=2-1]
	\arrow["{(\sigma(\chi_{a},\chi_{b}))^{\ast}}", from=1-4, to=2-4]
	\arrow["{\lambda^{\ast}_{\chi_{a}\chi_{b}}}"', from=2-1, to=2-4].
\end{tikzcd}\] one can conclude that 
\begin{align}
    \sigma(\chi_{a},\chi_{b})=\frac{c_{\chi_{a}}\lambda_{\chi_{b}\chi_{a}}}{\lambda_{\chi_{b}}}.
    %=\frac{\psi(\chi_{b})c_{\chi_{b}}\lambda^{'}_{\chi_{b}\chi_{a}}}{\lambda^{'}_{\chi_{a}}}.
\end{align}
From equation (\ref{15}) , we will get 
\begin{align}
    \sigma(\chi_{a},\chi_{b})=\frac{c_{\chi_{a}}\lambda_{\chi_{b}\chi_{a}}}{\lambda_{\chi_{b}}}=\frac{c_{\chi_{a}}c_{\chi_{b}}}{c_{\chi_{a}\chi_{b}}}.
\end{align} So, $\lambda_{\chi_{a}\chi_{b}}c_{\chi_{a}\chi_{b}}=\lambda_{\chi_{b}}c_{\chi_{b}}.$ If we assume $\chi_{b}=\chi_{000}$ then 
\begin{align}
    c_{\chi_{a}}\lambda_{\chi_{a}}=\lambda_{\chi_{000}}.
\end{align}
\end{proof}

\blmma
$\sigma(\chi_{a},\chi'_{b})=\frac{c_{\chi_{a}}e_{\chi'_{b}}}{e_{\chi_{a}\chi'_{b}}}$ and $ \sigma(\chi'_{a},\chi_{b})=\frac{c'_{\chi_{b}}e_{\chi'_{a}}}{e_{\chi'_{a}\chi_{b}}}$, where $\chi_a\in\{\chi_{i0k}\},\chi'_b\in\{\chi_{i1k}\}.$
\elmma
\begin{proof}
 \begin{align}
    \sigma(\chi_{a},\chi'_{b})=\frac{c_{\chi_{a}}e_{\chi'_{b}}}{e_{\chi_{a}\chi'_{b}}}
\end{align} follows from this commutative diagram
\[\begin{tikzcd}
	{\chi_{a}\otimes\chi^{'}_{b}\otimes\pi^{'}_{2}} &&& {\chi_{a}\otimes\pi_{2}} \\
	{\chi_{a}\chi^{'}_{b}\otimes\pi^{'}_{2}} &&& {\pi_{2}}
	\arrow["{(e_{\chi'_{b}})^{\ast}}", from=1-1, to=1-4]
	\arrow["{(\sigma(\chi_{a},\chi^{'}_{b}))^{\ast}}"', from=1-1, to=2-1]
	\arrow["{(c_{\chi_{a}})^{\ast}}", from=1-4, to=2-4]
	\arrow["{(e_{\chi_{a}\chi'_{b}})^{\ast}}"', from=2-1, to=2-4].
\end{tikzcd}\]
 One can easily observe that
 \begin{align}
      \sigma(\chi'_{a},\chi_{b})=\frac{c'_{\chi_{b}}e_{\chi'_{a}}}{e_{\chi'_{a}\chi_{b}}}
 \end{align}from
 \[\begin{tikzcd}
	{\chi'_a\ot\chi_b\ot\pi'_{2}} &&& {\chi'_a\ot \pi'_{2}} \\
	{\chi'_a\chi_b\ot\pi'_2} &&& {\pi_2}
	\arrow["{(c'_{\chi_2})^{\ast}}", from=1-1, to=1-4]
	\arrow["{(\sigma(\chi'_a,\chi_b)})^{\ast}"', from=1-1, to=2-1]
	\arrow["{(e_{\chi'_a})^{\ast}}", from=1-4, to=2-4]
	\arrow["{(e_{\chi'_a\chi_b})^{\ast}}", from=2-1, to=2-4].
\end{tikzcd}\]
   
\end{proof}

\blmma 
$e_{\chi'}^{2}$ is a constant.
\elmma
\begin{proof}
    For a fix $\chi\in\{\chi_{i0k}\},$ we denote $A_{\chi}:=(\Upsilon^{\pi\ot\pi}_{\chi} \ot \pi')\Upsilon^{\chi\ot\pi'}_{\pi'}$. Similarly, for a fixed $\chi'\in \{\chi_{i1k}\}$, denoted $B_{\chi'}:=(id_{\pi}\ot \Upsilon^{\pi\ot \pi'}_{\chi'})\Upsilon^{\pi\ot\chi'}_{\pi'}.$ One can easily check that $\{B_{\chi_{i1k}}:i,k=0,1\}$ is a basis of $Mor(\pi'_2,\pi_2\ot\pi_2\ot\pi'_2)$.

Assume that $\chi_{i}\in\{\chi_{m0n}:m,n\in\{0,1\}\}.$ We can write $A_{\chi_{i}}=\sum_{j} c_{ij}B_{\chi'_{j}},$ where $\chi'_{j}\in\{\chi_{m1n}\}, c_{ij}\in\bbc.$ From this, we will get
\begin{align}
 \widetilde{A}_{\chi_{i}}=\lambda_{\chi_{i}}c'_{\chi_{i}}A_{\chi_{i}}=\sum_{j}c_{ij}d_{\chi'_{j}}\eta_{\chi'_{j}}B_{\chi'_{j}}=\sum_{j} c_{ij}\widetilde{B}_{\chi'_{j}}.   
\end{align}
Now, we can conclude that
\begin{align}\label{30}
  \lambda_{\chi_{i}}c'_{\chi_{i}}= d_{\chi'_{j}}\eta_{\chi'_{j}}=a_0 ~for~all~\chi_{i}\in\{\chi_{mon}\},~\chi'_{j}\in \{\chi_{m1n}\},
\end{align} where $a_0$ is a constant.

Let $P_{\chi'}:=(\Upsilon^{\pi_2\ot\pi'_2}_{\chi'}\ot Id_{\pi'_2})\Upsilon^{\chi'\ot\pi'_2}_{\pi_2}$ and $Q_{\chi}:=(Id_{\pi_2}\ot \Upsilon^{\pi'_2\ot\pi'_2}_{\chi})\Upsilon^{\pi_2\ot\chi}_{\pi_2}.$ $\{Q_{\chi_{iok}}\}$ is a basis of $Mor(\pi_{2},\pi_2\ot\pi'_2\ot\pi'_2).$ 
Let us assume that $P_{\chi'_{i}}=\sum_{j} d_{ij}Q_{\chi_{j}}.$
Now,
\begin{align}
\widetilde{P_{\chi'_{i}}}=\eta_{\chi'_{i}}e_{\chi'_{i}}P_{\chi'_{i}}=\sum_{j}d_{ij}c_{\chi_{j}}\lambda'_{\chi_{j}}Q_{\chi_{j}}=\sum_{j} d_{ij}\widetilde{Q_{\chi_{j}}},
\end{align}
From which, one can easily observe that  for all $\chi'_{i}\in\{\chi_{m1n}\},\chi_{j}\in\{\chi_{m0n}\}$
\begin{align}\label{32}
 \eta_{\chi'_{i}}e_{\chi'_{i}}= c_{\chi_{j}}\lambda'_{\chi_{j}}=a'_{0},
\end{align} where $a'_{0}$ is a constant.

From equations (\ref{30}) and (\ref{32}), we can conclude that
\begin{align}
   d^{2}_{\chi'_i} =\frac{d_{\chi'_{j}}\eta_{\chi'_{j}}}{\eta_{\chi'_{i}}e_{\chi'_{i}}}=\frac{a_0}{a'_0}.
\end{align} Hence $e^{2}_{\chi'}$ is a constant.
\end{proof}

%\brmrk  From lemma (5.1.6), it follows that $\psi\equiv 1$ and now it is easy to check  $c'_{\chi}=c_{\chi}$ from lemma (5.1.4).

%\ermrk

%\brmrk From the equations (5.1.29) and (5.1.31), it follows that $a_{0}=a'_{0}.$
%\ermrk

 \blmma $\psi'\equiv 1$ on $\{\chi_{i0k}\}.$
 \elmma
 \begin{proof}
     From Lemma (\ref{lemma3.6}), it  folllows that $\psi'(\chi_{a})=1$, for all $\chi_{a}\in\{\chi_{i0k}\}.$
 \end{proof}

 Let $V_{\chi}:\bbc_{\chi}\to \bbc_{\chi}$, $V_{\chi'}:\bbc_{\chi'}\to \bbc_{\chi'},$  $V_{\pi_2}:\clh_{\pi_2}\to \clh_{\pi_2} $ and $V_{\pi'_2}:\clh_{\pi'_2}\to \clh_{\pi'_2}$ be  unitary linear maps defined by
 \begin{align}
    &V_{\chi}(1_{\bbc_{\chi}})=c_{\chi}1_{\bbc_{\chi}},~~ V_{\chi'}(1_{\bbc_{\chi'}})=\frac{a^{1/2}_0}{a'^{1/2}_0}e_{\chi'}1_{\bbc_{\chi'}},
   ~~V_{\pi_2}=a^{1/2}_{0}Id_{\clh_{\pi_2}} , ~~V_{\pi'_2}= a'^{1/2}_{0}Id_{\clh_{\pi'_2}}.
 \end{align}.
 \blmma $\sigma$ is equivalent to the trivial 2-cycle of the group $\{\chi_{ijk}\}.$\elmma
 \begin{proof}
 For any two $\chi_1,\chi_2\in\{\chi_{i0k}\},$ we have
 \begin{align*}
     \sigma(\chi_1,\chi_2)=&(V_{\chi_1}\otimes V_{\chi_2})V^{\ast}_{\chi_1\chi_2}\\
     =&c_{\chi_1}c_{\chi_2}(c_{\chi_1\chi_2})^{-1}.
 \end{align*}
 Assume that $\chi'_{1},\chi'_{2}\in\{\chi_{i1k}\}.$
 \begin{align*}
    \sigma(\chi_1,\chi'_2)=&(V_{\chi_1}\otimes V_{\chi'_2})V^{\ast}_{\chi_1\chi'_2}\\ 
    =&c_{\chi_1}\frac{a^{1/2}_0}{a'^{1/2}_0}e_{\chi'_2} 
     (e_{\chi_1\chi'_2})^{-1}\frac{a'^{1/2}_0}{a^{1/2}_0}\\
     =&\frac{c_{\chi_1}e_{\chi'_2}}{e_{\chi_1\chi'_{2}}},
 \end{align*}
 and 
 \begin{align*}
    \sigma(\chi'_1,\chi'_2)=&(V_{\chi'_1}\otimes V_{\chi'_2})V^{\ast}_{\chi'_1\chi'_2}\\ 
    =&(e_{\chi'_{1}}e_{\chi'_{2}}\frac{a_0}{a'_0})(c_{\chi'_1\chi'_2})^{-1}\\
    =&(\frac{e_{\chi'_1}}{c_{\chi'_1\chi'_2}})(e_{\chi'_2}\frac{a_0}{a'_0})\\
    =&\frac{e_{\chi'_1}}{c_{\chi'_1\chi'_2}}\frac{1}{e_{\chi'_2}}~~~~~~~~~~~~~~[~from~~e^{2}_{\chi'_2}=\frac{a'_0}{a_0}]\\
    =&\frac{e_{\chi'_1}d_{\chi'_2}}{c_{\chi'_1\chi'_2}}.     
 \end{align*}

 \begin{align*}
     \sigma(\chi'_1,\chi_2)=&(V_{\chi'_1}\otimes V_{\chi_2})V^{\ast}_{\chi'_1\chi_2}\\
     =& e_{\chi'_1}c_{\chi_2}(e_{\chi'_1\chi_2})^{-1}\\
     =&\frac{e_{\chi'_1}c_{\chi_2}}{e_{\chi'_1\chi_2}}\\
     =&\frac{e_{\chi'_1}c'_{\chi_2}}{e_{\chi'_1\chi_2}}. ~~~~[from~~c'_{\chi}=c_{\chi}]
 \end{align*}
 Hence, $\sigma$ is equivalent to the trivial 2-cycle.
 \end{proof}
% One can check that $\phi(\Upsilon^{a\ot b}_c)= (v_a \ot v_b)(\Upsilon^{a\ot b}_c) v^{\ast}_c,$ where $a,b,c\in \{\chi_{ijk},\pi_2,\pi'_2\}.$ Hence it is isomorphic to the identity tensor functor from the definition (\ref{isof}).

Without loss of generality assume that $\sigma$ is a trivial 2-cycle.
 
  Assume that $\phi(\Upsilon^{\chi\ot\pi_4}_{\pi_4})=\widetilde{\Upsilon}^{\chi\ot\pi_4}_{\pi_4}=\tau(\chi)\Upsilon^{\chi\ot\pi_4}_{\pi_4}.$
        \blmma $\tau$ is a character on the group $\{\chi_{ijk}\}.$ \elmma
        \begin{proof}
           For any $\chi_1,\chi_2\in\{\chi_{ijk}\},$ we will get 
           \begin{align}
               (\Upsilon^{\chi_1\ot\chi_2}_{\chi_{1}\chi_{2}} \ot Id_{\pi_4})\Upsilon^{\chi_{1}\chi_{2} \ot\pi_4}_{\pi_4}=(Id_{\chi_1}\ot \Upsilon^{\chi_2\ot\pi_4}_{\pi_4})\Upsilon^{\chi_1\ot\pi_4}_{\pi_4}.
           \end{align} 
           If we apply $\phi$ on both sides of the equation then \begin{align}
              \tau(\chi_1\chi_2)(\Upsilon^{\chi_1\ot\chi_2}_{\chi_{1}\chi_{2}} \ot Id_{\pi_4})\Upsilon^{\chi_{1}\chi_{2} \ot\pi_4}_{\pi_4} =\tau(\chi_1)\tau(\chi_2)(Id_{\chi_1}\ot \Upsilon^{\chi_2\ot\pi_4}_{\pi_4})\Upsilon^{\chi_1\ot\pi_4}_{\pi_4}.
           \end{align} Hence $\tau$ is a character on $\{\chi_{ijk}\}.$ 
        \end{proof}
Let us assume that $\widetilde{\Upsilon}^{\pi_4\ot\pi_4}_{\chi}:=m_{\chi}\Upsilon^{\pi_4\ot\pi_4}_{\chi}.$ Then the following identity holds,
           \begin{align*}
               (Id_{\chi_{1}} \ot \Upsilon^{\pi_4\ot\pi_4}_{\chi_2})\Upsilon^{\chi_{1}\ot\chi_{2}}_{\chi_1\chi_2}=(\Upsilon^{\chi_1\ot\pi_4}_{\pi_4}\ot Id_{\pi_4})\Upsilon^{\pi_4\ot\pi_4}_{\chi_1\chi_2}.
           \end{align*}
\blmma $m_{\chi_1}=\tau(\chi_1)m_{\chi_{000}}.$
\elmma
        \begin{proof}
        Observe that \begin{align}
               \phi((Id_{\chi_{1}} \ot \Upsilon^{\pi_4\ot\pi_4}_{\chi_2})\Upsilon^{\chi_{1}\ot\chi_{2}}_{\chi_1\chi_2}))=&(Id_{\chi_{1}} \ot \widetilde{\Upsilon}^{\pi_4\ot\pi_4}_{\chi_2})\widetilde{\Upsilon}^{\chi_{1}\ot\chi_{2}}_{\chi_1\chi_2}\\
               =& m_{\chi_2}(Id_{\chi_{1}} \ot \Upsilon^{\pi_4\ot\pi_4}_{\chi_2})\Upsilon^{\chi_{1}\ot\chi_{2}}_{\chi_1\chi_2}),
           \end{align}

            and
            \begin{align}
                \phi((\Upsilon^{\chi_1\ot\pi_4}_{\pi_4}\ot Id_{\pi_4})\Upsilon^{\pi_4\ot\pi_4}_{\chi_1\chi_2}))=&(\widetilde{\Upsilon}^{\chi_1\ot\pi_4}_{\pi_4} \ot Id_{\pi_4})\widetilde{\Upsilon}^{\pi_4\ot\pi_4}_{\chi_1\chi_2}\\
                =& \tau(\chi_1)m_{\chi_1\chi_2} ((\Upsilon^{\chi_1\ot\pi_4}_{\pi_4}\ot Id_{\pi_4})\Upsilon^{\pi_4\ot\pi_4}_{\chi_1\chi_2}).
            \end{align}
             We can conclude that $m_{\chi_2}=\tau(\chi_1)m_{\chi_1\chi_2}$ from equations (38) and (40). If we choose $\chi_{2}=\chi_{000}$ then $m_{\chi_1}=\tau(\chi_1)m_{\chi_{000}}.$ Assume that $\mu=m_{\chi_{000}}.$
        \end{proof} 

        \blmma $c_{\chi}=\tau(\chi)\psi(\chi)$ for a character $\psi$ on $\{\chi_{i0k}\}.$
        \elmma
\begin{proof}
From the  commutative diagram
    \[\begin{tikzcd}
	{\chi_{1}\otimes\chi_{2}\otimes\pi_4} & {\chi_{1}\otimes\pi_4} & {} \\
	{\chi_1\chi_2\otimes\pi_4} & {\pi_4}
	\arrow["{\tau^{\ast}(\chi_2)}", from=1-1, to=1-2]
	\arrow["{\sigma^{\ast}(\chi_1,\chi_2)}"', from=1-1, to=2-1]
	\arrow["{\tau^{*}(\chi_1)}", from=1-2, to=2-2]
	\arrow["{\tau^{*}(\chi_1\chi_2)}", from=2-1, to=2-2]
\end{tikzcd}\]
we have $\sigma(\chi_1,\chi_2)=\frac{\tau(\chi_1)\tau(\chi_2)}{\tau(\chi_1\chi_2)}=1$. From \eqref{15}, 
\begin{align*}
    \frac{c_{\chi_1}c_{\chi_2}}{c_{\chi_1\chi_2}}=\frac{\tau(\chi_1)\tau(\chi_2)}{\tau(\chi_1\chi_2)}.
\end{align*}
Hence there exists a character $\psi$ on $\{\chi_{i0k}\}$ such that $c_{\chi}=\tau(\chi)\psi(\chi)$ for $\chi\in \{\chi_{i0k}\}.$
\end{proof}

Let $\{T_1,T_2\}$ be a basis of $Mor(\pi_4, \pi_2\ot \pi_4).$ $T_1,T_2$ are defined by,
\begin{align*}
    T_{1}(x_1)=e_1\ot x_4, T_1(x_2)=e_2\ot x_3, T_1(x_3)=e_1\ot x_2, T_1(x_4)=e_2\ot x_1\\
     T_{2}(x_1)=e_2\ot x_2, T_2(x_2)=e_1\ot x_1, T_2(x_3)=e_2\ot x_4, T_2(x_4)=e_1\ot x_3.
\end{align*}
   Let $E^{\chi}_{i}:=(Id_{\chi}\ot T_i)\Upsilon^{\chi\ot\pi_4}_{\pi_4}\in Mor(\pi_4, \chi\ot\pi_2\ot\pi_4)$ and $F^{\chi}_{j}:=(\Upsilon^{\chi\ot\pi_2}_{\pi_2}\ot Id_{\pi_4})T_{j}\in Mor(\pi_4, \chi\ot\pi_2\ot\pi_4).$
    Now, one can easily check that the following equations hold,
    \begin{align}
         &E^{\chi_{000}}_1=F^{\chi_{000}}_1,  
        ~~~E^{\chi_{000}}_2=F^{\chi_{000}}_2,\\
         &E^{\chi_{001}}_1=F^{\chi_{001}}_2, ~~ E^{\chi_{001}}_2=F^{\chi_{001}}_1,\\
         &E^{\chi_{101}}_1=-F^{\chi_{101}}_2,~~ E^{\chi_{101}}_2=-F^{\chi_{101}}_1,\\
           &E^{\chi_{100}}_1=F^{\chi_{100}}_1,~~ E^{\chi_{100}}_2=-F^{\chi_{100}}_2.
        \end{align}

        Let us assume that $\widetilde{T_1}=\sum^{2}_{i=1} \omega_{1k} T_{k}, \widetilde{T_2}=\sum^{2}_{i=1} \omega_{2k} T_{k},$ where $\omega_{ij}\in \bbc, i,j\in\{1,2\}$ and $E^{\chi}_{i}=\sum_{k} a^{\chi}_{ik} F^{\chi}_{k}.$
        \blmma $\tau(\chi)\omega a^{\chi}=c_{\chi}a^{\chi}\omega.$ \elmma
        \begin{proof}
        \begin{align}
            \widetilde{F^{\chi}_{j}}=&c_{\chi}(\Upsilon^{\chi\ot\pi_2}_{\pi_2} \ot Id_{\pi_4})\widetilde{T_j}\\
            =&c_{\chi}\sum_{k}\omega_{jk} (\Upsilon^{\chi\ot\pi_2}_{\pi_2} \ot Id_{\pi_4}) T_{k}\\
            =&c_{\chi}\sum_{k}\omega_{jk}F^{\chi}_{k}.
        \end{align}
 Now,
        \begin{align}
            \phi(E^{\chi}_{i})=\widetilde{E^{\chi}_{i}}=&\sum_{k} a^{\chi}_{ik}\widetilde{F^{\chi}_{k}}\\
            =&c_{\chi}\sum_{k,p} a^{\chi}_{ik}\omega_{kp} F^{\chi}_p
        \end{align}
       and
        \begin{align}
            \widetilde{E^{\chi}_{j}}=& (Id_{\chi} \otimes \widetilde{T_j})\tau(\chi)\Upsilon^{\chi\ot \pi_4}_{\pi_4}\\
            =&\sum_{k} \omega_{jk}(Id_{\chi}\ot T_k)\tau(\chi)\Upsilon^{\chi\ot\pi_4}_{\pi_4}\\
            =&\sum_{k}\omega_{jk}\tau(\chi) E^{\chi}_{k}\\
            = &\sum_{k,l}\omega_{jk}\tau(\chi)a^{\chi}_{kl} F^{\chi}_{l}\\
            =& \tau(\chi)\sum_{k,l} \omega_{jk}a^{\chi}_{kl}F^{\chi}_{l}.
        \end{align}
         From (49) and (54), we will get 
         \begin{align}\label{55}
            \tau(\chi)\omega a^{\chi}=c_{\chi}a^{\chi}\omega,
         \end{align} where $\omega=(\omega_{ij}), a^{\chi}=(a_{ij})$ are 2$\times$ 2 matrices.
\end{proof}
        From (41),(42),(43) and (44), 
        \begin{align}
        a^{\chi_{000}}=\begin{pmatrix}
            1 & 0\\
            0 & 1
        \end{pmatrix},
            a^{\chi_{001}}=\begin{pmatrix}
                0 & 1\\
                1 &0
            \end{pmatrix},
           a^{\chi_{101}} =\begin{pmatrix}
               0& -1\\
               -1 & 0
           \end{pmatrix},
           a^{\chi_{100}}=\begin{pmatrix}
               1 & 0\\
               0 & -1
           \end{pmatrix}.
        \end{align}

        \blmma There are only two choices of $\psi$, which are $\psi\equiv 1$ and $\psi(\chi_{100})=1, \psi(\chi_{001})=-1$ and the followings hold
        \begin{enumerate}
            \item [1)] If $\psi\equiv 1$ then $\omega=\lambda Id,$ where $\lambda\in \bbc-\{0\}.$
            \item[2)]  If $\psi(\chi_{100})=1,\psi(\chi_{001})=-1$ then $\omega=diag( \lambda,-\lambda)$ for a $\lambda\in \bbc-\{0\}.$
           % \item[3)] If $\tau(\chi_{100})= -1,\tau(\chi_{001})=1$ then \begin{align*}
               % \omega=\begin{pmatrix}
                %0 & \lambda\\
                %\lambda & 0
            %\end{pmatrix}.
            %\end{align*}
            %\item[4)] If $\tau(\chi_{100})= -1,\tau(\chi_{001})= -1$ then \begin{align*}
                %\omega=\begin{pmatrix}
               % 0 & \lambda\\
                %-\lambda & 0
            %\end{pmatrix}.
            %\end{align*}
        \end{enumerate}
        \elmma
        \begin{proof}
From the equation \eqref{55}, one can observe that 
\begin{align*}
    \tau(\chi)wa^{\chi}=&c_{\chi}a^{\chi}w\\
    =&\tau(\chi)\psi(\chi)a^{\chi}w.
\end{align*}  So, it is easily follows $\psi(\chi)wa^{\chi}=a^{\chi}w.$
\begin{enumerate}
    \item [1)] For $\psi\equiv1,$ $w$ commutes with the matrices $a^{\chi_{000}},a^{\chi_{101}},a^{\chi_{101}},a^{\chi_{100}}.$ From $wa^{\chi_{100}}=a^{\chi_{100}}w$ we can conclude $w=diag(a,b)$ and from $wa^{\chi_{001}}=a^{\chi_{001}}w$ one can check that $a=b.$ Hence $w=diag(\lambda,\lambda)$ for a non-zero complex number $\lambda.$
    \item[2)] If $\psi(\chi_{100})=1,\psi(\chi_{001})=-1$ then $w=diag(a,b)$ from previous observation. From $(-1)wa^{\chi_{001}}=a^{\chi_{001}}w$ one can prove that $b=-a$ and it satisfies $(-1)diag(a, -a)a^{\chi_{101}}=a^{\chi_{101}}diag(a,-a).$ So $w=diag(\lambda, -\lambda)$ for a $\lambda\in\bbc-\{0\}.$
    \item[3)] If $\psi(\chi_{100})=-1,\psi(\chi_{001})=1$ then $w=\begin{pmatrix}
        0 & b\\
        c &  0
    \end{pmatrix}$ from the equation $\psi(\chi_{100})wa^{\chi_{100}}=a^{\chi_{100}}w.$ It is easy to observe $w=\begin{pmatrix}
        0 & b\\
        b & 0
    \end{pmatrix}$ from equation $\phi(\chi_{001})\begin{pmatrix}
        0 & b\\
        c &  0
    \end{pmatrix}a^{\chi_{001}}=a^{\chi_{001}}\begin{pmatrix}
        0 & b\\
        c & 0
    \end{pmatrix}.$ Now, we can conclude $b=0$ from equation $\psi(\chi_{101})\begin{pmatrix}
        0 & b\\
        b & 0
    \end{pmatrix}a^{\chi_{101}}=a^{\chi_{101}}\begin{pmatrix}
        0 & b\\
        b & 0
    \end{pmatrix}.$ So there is no such matrix $w$ for $\psi(\chi_{100})=-1,\psi(\chi_{001})=1.$
    \item[4)] Similarly we can prove that there is no such matrix $w$ for $\psi(\chi_{100})=-1,\psi(\chi_{001})=-1.$
\end{enumerate}

        \end{proof}

      Let  $T'_1,T'_2\in Mor(\pi_4, \pi'_2\ot \pi_4)$ be  defined by,
\begin{align*}
    T'_{1}(x_1)=f_1\ot x_4, T'_1(x_2)=f_2\ot x_3, T'_1(x_3)=-f_1\ot x_2, T_1(x_4)=-f_2\ot x_1\\
     T'_{2}(x_1)=f_2\ot x_2, T'_2(x_2)=f_1\ot x_1, T'_2(x_3)=-f_2\ot x_4, T'_2(x_4)=-f_1\ot x_3.
\end{align*} Then $\{T'_1,T'_2\}$ is a basis of $Mor(\pi_4, \pi'_2\ot \pi_4).$
 Let $E^{\chi'}_{i}:=(Id_{\chi'}\ot T_i)\Upsilon^{\chi'\ot\pi_4}_{\pi_4}\in Mor(\pi_4, \chi\ot\pi_2\ot\pi_4)$ and $F^{\chi'}_{j}:=(\Upsilon^{\chi'\ot\pi_2}_{\pi'_2}\ot Id_{\pi_4})T'_{j}\in Mor(\pi_4, \chi\ot\pi_2\ot\pi_4).$
    Now, one can easily check that the following equations hold,
    \begin{align}
         &E^{\chi_{010}}_1=F^{\chi_{010}}_1,  
        ~~~E^{\chi_{010}}_2=F^{\chi_{010}}_2,\\
         &E^{\chi_{011}}_1=-F^{\chi_{011}}_2, ~~ E^{\chi_{011}}_2=-F^{\chi_{011}}_1,\\
         &E^{\chi_{111}}_1=-F^{\chi_{111}}_2,~~ E^{\chi_{111}}_2=F^{\chi_{111}}_1,\\
           &E^{\chi_{110}}_1=F^{\chi_{110}}_1,~~ E^{\chi_{110}}_2=-F^{\chi_{110}}_2.
        \end{align}
Assume that $\widetilde{T'_1}=\sum^{2}_{i=1} \omega'_{1k} T'_{k}, \widetilde{T'_2}=\sum^{2}_{i=1} \omega'_{2k} T'_{k},$ where $\omega'_{ij}\in \bbc, i,j\in\{1,2\}$ and $E^{\chi'}_{i}=\sum_{k} a^{\chi'}_{ik} F^{\chi'}_{k}.$
\blmma 
$\tau(\chi')\omega a^{\chi'}= d_{\chi'}a^{\chi'}\omega'.$
\elmma
\begin{proof}
    \begin{align}
           \phi(F^{\chi'}_{j})= \widetilde{F^{\chi'}_{j}}=&d_{\chi'}(\Upsilon^{\chi'\ot\pi_2}_{\pi'_2} \ot Id_{\pi_4})\widetilde{T'_j}\\
            =&d_{\chi'}\sum_{k}\omega'_{jk} (\Upsilon^{\chi'\ot\pi_2}_{\pi'_2} \ot Id_{\pi_4}) T'_{k}\\
            =&d_{\chi'}\sum_{k}\omega'_{jk}F^{\chi'}_{k}.
        \end{align}
        
        \begin{align}
            \phi(E^{\chi'}_{i})=\widetilde{E^{\chi'}_{i}}=&\sum_{k} a^{\chi'}_{ik}\widetilde{F^{\chi'}_{k}}\\
            =& d_{\chi'}\sum_{k,p} a^{\chi'}_{ik}\omega'_{kp} F^{\chi'}_p,
        \end{align}
       and
        \begin{align}
            \widetilde{E^{\chi'}_{j}}=& (Id_{\chi'} \otimes \widetilde{T_j})\tau(\chi')\Upsilon^{\chi'\ot \pi_4}_{\pi_4}\\
            =&\sum_{k} \omega_{jk}(Id_{\chi'}\ot T_k)\tau(\chi')\Upsilon^{\chi'\ot\pi_4}_{\pi_4}\\
            =&\sum_{k}\omega_{jk}\tau(\chi') E^{\chi'}_{k}\\
            = &\sum_{k,l}\omega_{jk}\tau(\chi')a^{\chi'}_{kl} F^{\chi'}_{l}\\
            =& \tau(\chi')\sum_{k,l} \omega_{jk}a^{\chi'}_{kl}F^{\chi'}_{l}.
        \end{align}
         From equations (65) and (70), we will get 
         \begin{align}\label{77}
            \tau(\chi')\omega a^{\chi'}=d_{\chi'}a^{\chi'}\omega',
         \end{align} where $\omega'=(\omega'_{ij}), a^{\chi'}=(a'_{ij})$ are 2$\times$ 2 matrices.
         
\end{proof}
 From equations (57),(58),(59) and (60), One can conclude that 
       \begin{align}\label{78}
        a^{\chi_{010}}=\begin{pmatrix}
            1 & 0\\
            0 & 1
        \end{pmatrix},
            a^{\chi_{111}}=\begin{pmatrix}
                0 & 1\\
                -1 &0
            \end{pmatrix},
           a^{\chi_{011}} =\begin{pmatrix}
               0& -1\\
               -1 & 0
           \end{pmatrix},
           a^{\chi_{110}}=\begin{pmatrix}
               1 & 0\\
               0 & -1
           \end{pmatrix}.
        \end{align}
   If $\psi\equiv 1$ then from equation \eqref{77} one can observe that 
   \begin{align*}
       w'=&\tau(\chi')e_{\chi'}(a^{\chi'})^{-1}w a^{\chi'}\\
       =&\tau(\chi')e_{\chi'}w\\
       =&\tau(\chi')e_{\chi'}\lambda Id.
   \end{align*} So there exists a constant $k_{0}$ such that $\tau(\chi')e_{\chi'}=k_{0}$ for all $\chi'\in\{\chi_{ijk}\}$ and also $k^{2}_{0}=\frac{a'_{0}}{a_0}.$

If $\psi(\chi_{100})=1,\psi(\chi_{001})=-1$ then \begin{align*}
    w'=&\tau(\chi')e_{\chi'}(a^{\chi'})^{-1}w a^{\chi'}.\\
\end{align*}
From the matrices \eqref{78}, it is a routine check that
\begin{enumerate}
    \item [1)] For $\chi'=\chi_{010}$, $w'=\tau(\chi_{010})e_{\chi_{010}}\begin{pmatrix}
        \lambda & 0\\
        0    &  -\lambda
    \end{pmatrix},$
    \item[2)] For $\chi'=\chi_{110}, w'=e_{\chi_{110}}\tau(\chi_{110})\begin{pmatrix}
        \lambda & 0\\
        0 & -\lambda.
    \end{pmatrix},$
    \item[3)] For $\chi'=\chi_{011}, w'= -e_{\chi_{011}}\tau(\chi_{011})\begin{pmatrix}
        \lambda & 0\\
        0 & -\lambda.
    \end{pmatrix},$
    \item[4)] For $\chi'=\chi_{111}, w'= -e_{\chi_{111 }}\tau(\chi_{111})\begin{pmatrix}
        \lambda & 0\\
        0 & -\lambda.
    \end{pmatrix}.$
\end{enumerate}
So, $\tau(\chi_{010})e_{\chi_{010}}=\tau(\chi_{110})e_{\chi_{110}}=-\tau(\chi_{011})e_{\chi_{011}}=-\tau(\chi_{111})e_{\chi_{111}}.$
       % \begin{enumerate}
          %  \item [1)] If $\psi\equiv 1$ then from equation (6.1.70) one can conclude that $\omega a^{\chi'}=a^{\chi'}\omega^{'}.$ Hence, $\omega'=\omega= diag  (\lambda, \lambda).$
            %\item[2)] If $\tau(\chi_{010})=1, \tau(\chi_{100})=1,\tau(\chi_{001})=-1$ then $\omega'=\omega=diag( \lambda,-\lambda).$
            %\item[3)] If $\tau(\chi_{010})= -1, \tau(\chi_{100})=1,\tau(\chi_{001})=-1$ then $\omega'= diag(-\lambda, \lambda) .$
            %\item[4)]  If $\tau(\chi_{010})= -1, \tau(\chi_{100})=1,\tau(\chi_{001})=1$ then $\omega'= -\omega= diag( -\lambda, -\lambda).$
            %\item[5)]  If $\tau(\chi_{010})=-1, \tau(\chi_{100})=-1,\tau(\chi_{001})=1$ then \begin{align*}
                %\omega'=\begin{pmatrix}
                 %   0 & -\lambda\\
                  %  -\lambda & 0
                %\end{pmatrix}.
           % \end{align*}
            %\item[6)]  If $\tau(\chi_{010})=1, \tau(\chi_{100})=-1,\tau(\chi_{001})=1$ then \begin{align*}
              %  \omega'=\begin{pmatrix}
            % $$       0 & \lambda\\
             %       \lambda & 0
              %  \end{pmatrix}.
            %\end{align*}
           % \item[7)]  If $\tau(\chi_{010})=1, \tau(\chi_{100})=-1,\tau(\chi_{001})= -1$ then \begin{align*}
              %  \omega'=\begin{pmatrix}
              %      0 & \lambda\\
              %      -\lambda & 0
               % \end{pmatrix}.
           % \end{align*}
            %\item[8)]  If $\tau(\chi_{010})=-1, \tau(\chi_{100})=-1,\tau(\chi_{001})=-1$ then \begin{align*}
              %  \omega'=\begin{pmatrix}
               %     0 & -\lambda\\
                %    \lambda & 0
               % \end{pmatrix}.
           % \end{align*}
       % \end{enumerate}
        \blmma $\lambda^2=a_{0}.$
        \elmma
         \begin{proof} The space $Mor(\pi_4, ~ \pi_2\ot \pi_2\ot \pi_4)$ is 4-dimensional, observe that $C_{\chi_{i0k}}=(\Upsilon^{\pi_2\ot\pi_2}_{\chi_{i0k}}\ot Id_{\pi_4})\Upsilon^{\chi_{i0k}\ot\pi_4}_{\pi_4}.$ $\{C_{\chi_{i0k}}\}$ is a basis of $Mor(\pi_4, ~ \pi_2\ot \pi_2\ot \pi_4),$ where
         \begin{align}
             C_{\chi_{000}}=(Id_{\pi_2}\ot T_1)T_1 +(Id_{\pi_2}\ot T_2)T_2,\\
              C_{\chi_{100}}=(Id_{\pi_2}\ot T_1)T_1 -(Id_{\pi_2}\ot T_2)T_2,\\
               C_{\chi_{001}}=(Id_{\pi_2}\ot T_1)T_2 +(Id_{\pi_2}\ot T_2)T_1,\\
                C_{\chi_{101}}=(Id_{\pi_2}\ot T_1)T_2 -(Id_{\pi_2}\ot T_2)T_1.
         \end{align}
         One can observe that \begin{align*}
             \widetilde{C_{\chi_{i0k}}}=\lambda_{\chi_{i0k}}\tau(\chi_{i0k})C_{\chi_{i0k}}=&a_{o}c_{\chi_{i0k}}\tau(\chi_{i0k})C_{\chi_{i0k}}\\
             =& a_{0}\tau(\chi_{i0k})\psi(\chi_{iok})\tau(\chi_{i0k})C_{\chi_{i0k}}\\
             =& a_{0}\psi(\chi_{i0k})C_{\chi_{i0k}}.
         \end{align*}
          % From equation $(5.1.78)$, one can observe that
          % \begin{align}
             %&\widetilde{C_{\chi_{000}}}=\tau(\chi_{000})C_{\chi_{000}}=C_{\chi_{000}}=(Id_{\pi_2}\ot \widetilde{T_1})\widetilde{T_1} +(Id_{\pi_2}\ot \widetilde{T_2})\widetilde{T_2}.
          % \end{align}
           From equations (73),(74), (75) and (76) it follows that
           \begin{enumerate}
               \item [i)] For $\psi\equiv 1,$
               %\widetilde{T_1}=\lambda T_1~~ and~~ \widetilde{T_2} =\lambda T_2.$ 
               We will get $\lambda^{2}= a_0.$
               \item[ii)] For $\psi(\chi_{100})=1,~\psi(\chi_{001})= -1,$  we know $\widetilde{T_1}=\lambda T_1,~\widetilde{T_2}=-\lambda T_2$. Now, it is a straightforward computation to prove  $\lambda^{2}=a_0.$
                %%%\item[iii)] For $\tau(\chi_{100})=-1,~\tau(\chi_{001})= 1,$$\widetilde{T_1}=\lambda T_2,~\widetilde{T_2}=\lambda T_1$. If we will put this in equation (5.1.82), we will get $\lambda^{2}=1.$
               % \item[iv)] For $\tau(\chi_{100})=-1,~\tau(\chi_{001})= -1,$  we know $\widetilde{T_1}=-\lambda T_2,~\widetilde{T_2}=\lambda T_1$. From equation (5.1.82), we can observe that $\lambda^{2}=1.$
           \end{enumerate}

           %Similarly, one can conclude that $\lambda^{2}=1$ from equations (5.1.79),(5.1.80) and (5.1.81).
           \end{proof}

             Let $S_{1},S_{2}$ be two linear maps defined by, 
             \begin{align*}
                 S_{1}(e_{1})=x_1\ot x_1+ x_3\ot x_3,~~ S_{1}(e_{2})=x_2\ot x_2 + x_4\ot x_4,\\
                 S_2(e_1)=x_2\ot x_4 +x_4\ot x_2,~~ S_{2}(e_2)= x_1\ot x_3 +x_3\ot x_1.
             \end{align*} $\{S_1,S_2\}$ is a basis of the vector space $Mor(\pi_2, \pi_4\ot \pi_4).$
             
             Let us define  \begin{align}
                 &G^{\chi}_{i}:=(\Upsilon ^{\chi\ot\pi_4}_{\pi_4} \ot Id_{\pi_4})S_{i},\\
                 &H^{\chi}_{i}:=(Id_{\chi} \ot S_{i})\Upsilon^{\chi\ot\pi_2}_{\pi_2},
             \end{align} where $G^{\chi}_{i},H^{\chi}_{i}\in Mor(\pi_2, \chi \ot \pi_4\ot\pi_4).$

             It is straightforward to verify that 
             \begin{align}
                 &G^{\chi_{000}}_{1}=H^{\chi_{000}}_1,~~G^{\chi_{000}}_2=H^{\chi_{000}}_2 ,\\
                 &G^{\chi_{001}}_{1}=H^{\chi_{001}}_{2},~~G^{\chi_{001}}_{2}=H^{\chi_{001}}_1,\\
                 &G^{\chi_{100}}_1=H^{\chi_{100}}_1,~~G^{\chi_{100}}_2= -H^{\chi_{100}}_2,\\
                 &G^{\chi_{101}}_1=-H^{\chi_{101}}_2,~~G^{\chi_{101}}_{2}= H^{\chi_{101}}.
                 \end{align}
            Let $\phi(S_l)=\widetilde{S_l}=\sum^{2}_{m=1} \theta_{lm}S_m$, where $\theta_{lm}\in \bbc,l=1~ to~2$ and also assume that $G^{\chi_{i0k}}_l=\sum^{2}_{y=1} n^{\chi_{i0k}}_{ly}H^{\chi_{i0k}}_{y},$ where $n^{\chi_{i0k}}_{ly}\in\bbc.$
\blmma \begin{align*}
    \tau(\chi_{i0k}) \theta n^{\chi_{i0k}}=c_{\chi_{i0k}}n^{\chi_{i0k}}\theta,
            \end{align*} where $\theta=(\theta_{lm}), n^{\chi_{i0k}}=(n^{\chi_{i0k}}_{ly})$ are  $2\times 2 $ matrices.
            \elmma

       \begin{proof}
            We have, \begin{align}
                \widetilde{G^{\chi_{i0k}}_{l}}=&\sum^{2}_{y=1} n^{\chi_{i0k}}_{ly}\widetilde{H^{\chi_{i0k}}_{y}}\\
                =& c_{\chi_{i0k}}\sum^{2}_{y=1} n^{\chi_{i0k}}_{ly}(Id_{\chi_{i0k}}\otimes \widetilde{S_y})\Upsilon^{\chi_{i0k}\ot\pi_2}_{\pi_2}\\
                =& c_{\chi_{i0k}}\sum^{2}_{y,z=1} n^{\chi_{i0k}}_{ly}\theta_{yz}(Id_{\chi_{i0k}}\otimes S_z)\Upsilon^{\chi_{i0k}\ot\pi_2}_{\pi_2}\\
                =& c_{\chi_{i0k}}\sum_{y,z}n^{\chi_{i0k}}_{ly}\theta_{yz}H^{\chi_{i0k}}_z,
            \end{align}
            
            and also \begin{align}
                \widetilde{G^{\chi_{i0k}}_{l}}=&(\widetilde{\Upsilon^{\chi_{i0k}\ot \pi_4}_{\pi_4}} \ot Id_{\pi_4})\widetilde{S_l}\\
                =&\tau(\chi_{i0k})(\Upsilon^{\chi_{i0k}\ot \pi_4}_{\pi_4} \ot Id_{\pi_4})\sum^{2}_{m=1} \theta_{lm}S_m\\
                =&\tau(\chi_{i0k})\sum_{m} \theta_{lm}(\Upsilon^{\chi_{i0k}\ot \pi_4}_{\pi_4} \ot Id_{\pi_4})S_m\\
                =&\tau(\chi_{i0k})\sum_{m} \theta_{lm} G^{\chi_{i0k}}_{m}\\
                =&\tau(\chi_{i0k})\sum_{m,z} \theta_{lm} n^{\chi_{i0k}}_{mz} H^{\chi_{i0k}}_{z}.
            \end{align}

            After comparing equations (86) and (91), one can conclude that \begin{align}\label{98}
            \tau(\chi_{i0k}) \theta n^{\chi_{i0k}}=c_{\chi_{i0k}}n^{\chi_{i0k}}\theta,
            \end{align} where $\theta=(\theta_{lm}), n^{\chi_{i0k}}=(n^{\chi_{i0k}}_{ly})$ are  $2\times 2 $ matrices.
       \end{proof}    
       \brmrk From the equation \eqref{98}, it is easy to conclude that $\psi(\chi_{i0k})\theta n^{\chi_{i0k}}=n^{\chi_{i0k}}\theta.$
       \ermrk
     We already know that 
     \begin{align}
        n^{\chi_{000}}=\begin{pmatrix}
            1 & 0\\
            0 & 1
        \end{pmatrix},
            n^{\chi_{001}}=\begin{pmatrix}
                0 & 1\\
                1 &0
            \end{pmatrix},
           n^{\chi_{101}} =\begin{pmatrix}
               0& 1\\
               -1 & 0
           \end{pmatrix},
           a^{\chi_{100}}=\begin{pmatrix}
               1 & 0\\
               0 & -1
           \end{pmatrix}.
        \end{align}
        
            \begin{enumerate}
                \item [1)] If $\psi\equiv 1$ then $\theta$ commutes with $n^{\chi_{100}}, n^{\chi_{101}}$. So, $\theta= diag (a, b)$ for two complex numbers $a,b.$ $\theta $ also commutes with $n^{\chi_{001}}$, from that one can easily reduce that $\theta= diag (\lambda_1, \lambda_1)$ for a non-zero complex number $\lambda_1.$
                \item[2)] If $\psi(\chi_{100})=1, \psi(\chi_{001})=-1$ then $\theta= diag (\lambda_1, -\lambda_1)$ where $\lambda_1\in \bbc-\{0\}.$
            \end{enumerate}

             Similarly, there is a basis $\{S'_1,S'_2\}$ of the vector space $Mor(\pi'_2, \pi_4\ot \pi_4).$

             Here we introduce another set of notations:
             \begin{align}
                 &G^{\chi_{i1k}}_{l'}:=(\Upsilon^{\chi_{i1k}\ot \pi_4}_{\pi_4} \ot Id_{\pi_4})S'_l\\
                 & H^{\chi_{i1k}}_{l}:=(Id_{\chi_{i1k}} \ot S_{l})\Upsilon^{\chi_{i1k}\ot\pi_2}_{\pi'_2},
             \end{align} where $G^{\chi_{i1k}}_{l'},H^{\chi_{i1k}}_{l}\in Mor(\pi'_2, \chi_{i1k}\ot \pi_4\ot \pi_4).$
             
     Assume that \begin{align}
         &\phi(S'_{l})=\widetilde{S'_{l}}=\sum_{m'}\theta'_{l'm'}S'_{m},\\
         & G^{\chi_{i1k}}_{m'}=\sum_{y} n'^{\chi_{i1k}}_{m'y} H^{\chi_{i1k}}_{y},
     \end{align}  where $\theta'_{l'm'},n'^{\chi_{i1k}}_{m'y}\in\bbc.$
        \begin{align*}
           \phi(G^{\chi_{i1k}}_{l'})=\widetilde{G^{\chi_{i1k}}_{l'}}=\tau(\chi_{i1k})((\Upsilon^{\chi_{i1k}\ot \pi_4}_{\pi_4} \ot \pi_4)\widetilde{S'_{l}} ).
        \end{align*}

         Similarly, we can observe that 
         \begin{align}
             \tau(\chi_{i1k})\theta'n'^{\chi_{i1k}}=d_{\chi_{i1k}}n'^{\chi_{i1k}}\theta.
         \end{align}
        \begin{enumerate}
         
            \item[1)] If $\psi\equiv 1$ then $\theta'=\tau(\chi_{i1k})d_{\chi_{i1k}}\lambda_1 Id,$ where $\tau(\chi_{i1k})d_{\chi_{i1k}}=\frac{1}{k_0}.$
            %$\theta'n'^{\chi_{i1k}} = n'^{\chi_{i1k}}\theta.$ From this, one can easily observe that $\theta'=\theta= diag (\lambda_1, \lambda_1).$
            \item[2)] If $\psi(\chi_{100})=1,\psi(\chi_{001})= -1,$
            %\tau(\chi_{010})= -1$
            then \item [i)] For $\chi'=\chi_{010}$, $\theta'=\tau(\chi_{010})d_{\chi_{010}}\begin{pmatrix}
        \lambda_1 & 0\\
        0    &  -\lambda_1
    \end{pmatrix},$
    \item[ii)] For $\chi'=\chi_{110}, \theta'=d_{\chi_{110}}\tau(\chi_{110})\begin{pmatrix}
        \lambda_1 & 0\\
        0 & -\lambda_1.
    \end{pmatrix},$
    \item[iii)] For $\chi'=\chi_{011}, \theta'= -d_{\chi_{011}}\tau(\chi_{011})\begin{pmatrix}
        \lambda_1 & 0\\
        0 & -\lambda_1.
    \end{pmatrix},$
    \item[iv)] For $\chi'=\chi_{111}, \theta'= -d_{\chi_{111 }}\tau(\chi_{111})\begin{pmatrix}
        \lambda_1 & 0\\
        0 & -\lambda_1.
    \end{pmatrix}.$

Hence $\tau(\chi_{010})e_{\chi_{010}}=\tau(\chi_{110})e_{\chi_{110}}=-\tau(\chi_{011})e_{\chi_{011}}=-\tau(\chi_{111})e_{\chi_{111}}.$
\end{enumerate}
            %$\theta'= -\theta= diag( -\lambda_1, -\lambda_1).$
            %\item[3)] If $\tau(\chi_{100})=1,\tau(\chi_{001})= -1, \tau(\chi_{010})= 1$ then $\theta'= \theta= diag( \lambda_1, -\lambda_1).$
           % \item[4)] If $\tau(\chi_{100})=1,\tau(\chi_{001})= -1, \tau(\chi_{010})= -1$ then $\theta'= -\theta= diag( -\lambda_1, \lambda_1).$

        Let us assume that \begin{align*}
            K_{i0k}:=(\Upsilon^{\pi_4\ot\pi_4}_{\chi_{iok}} \ot Id_{\pi_4})\Upsilon^{\pi_2\ot\chi_{i0k}}_{\pi_2}\in Mor(\pi_2, \pi_2\ot\pi_4\ot \pi_4).
        \end{align*}

    We already know that $\{S_i:i=1,2\}$ is a basis of $Mor(\pi_2, \pi_4\ot\pi_4)$ and $\{T_i: i=1,2\}$ is a bais of $Mor(\pi_4, \pi_2 \ot \pi_4).$
    Let \begin{align*}
        D_{ij}=(T_i \otimes Id_{\pi_4})S_j\in Mor(\pi_2, \pi_2 \ot\pi_4\ot\pi_4).
    \end{align*}
    Then the following relations hold,
    \begin{align}
        &D_{11}= (K_{000}-K_{100}),\\
        &D_{12}= (K_{001}-K_{101}),\\
        &D_{21}= (K_{001}+K_{101}),\\
        &D_{22}= (K_{000}+K_{100}).
    \end{align}
\blmma $\lambda_1\lambda=\mu.$
\elmma
  \begin{proof}   \begin{align}
        \phi(D_{11})=\widetilde{D_{11}}=&(\widetilde{T_1} \ot Id_{\pi_4})\widetilde{S_1}\\
        =&(\mu \tau(\chi_{000})c_{\chi_{000}} K_{000}- \mu \tau(\chi_{100})c_{\chi_{100}} K_{100}),
    \end{align} from which we can conclude that
    \begin{align}
        \widetilde{S_1}= (\widetilde{T^{\ast}_1} \ot Id_{\pi_4})\mu(\psi(\chi_{000})K_{000}-\psi(\chi_{100})K_{001}).
    \end{align}
     Similarly we can observe that
     \begin{align}
         &\widetilde{S_1}= (\widetilde{T^{\ast}_2} \ot Id_{\pi_4})\mu(\psi(\chi_{001})K_{001}+\psi(\chi_{101})K_{101}),\\
         &\widetilde{S_2}= (\widetilde{T^{\ast}_1} \ot Id_{\pi_4})\mu(\psi(\chi_{001})K_{001}-\psi(\chi_{101})K_{101}),\\
         &\widetilde{S_2}= (\widetilde{T^{\ast}_2} \ot Id_{\pi_4})\mu(\psi(\chi_{000})K_{000}+\psi(\chi_{100})K_{100}).
         \end{align}
         \begin{enumerate}
             \item [1)] If $\psi\equiv1$ then \begin{align*}
                 \widetilde{S_1}=&(\widetilde{T^{\ast}_1} \ot Id_{\pi_4})\mu(\tau(\chi_{000})K_{000}-\tau(\chi_{100})K_{001})\\
                 =&\frac{1}{\lambda}( T^{\ast}_1\ot Id_{\pi_4})\mu (K_{000}-K_{100})=\frac{\mu}{\lambda} S_1.
             \end{align*}
             Similarly, we can prove that $\widetilde{S_{2}}=\frac{\mu}{\lambda} S_2.$
             %\item[2)] If $\tau(\chi_{010})=-1, \tau(\chi_{100})=1,\tau(\chi_{001})=1$ then $\widetilde{S_1}=\lambda\mu S_1$ and $\widetilde{S_2}= \lambda\mu S_2.$
             \item[2)]If $ \psi(\chi_{100})=1,\psi(\chi_{001})=-1$ then $\widetilde{S_1}=\frac{\mu}{\lambda} S_1$ and $\widetilde{S_2}= -\frac{\mu}{\lambda} S_2.$
            % \item[4)] If $\tau(\chi_{010})=-1, \tau(\chi_{100})=1,\tau(\chi_{001})=-1$ then $\widetilde{S_1}=\lambda\mu S_1$ and $\widetilde{S_2}= -\lambda\mu S_2.$
         \end{enumerate}
    From this, we can easily conclude that $\lambda_1\lambda= \mu.$
    \end{proof}
  %  \bthm If $\tau(\chi_{100})=1,\tau(\chi_{001})=-1 ,\tau(\chi_{010})= 1$ then it's correspondence to a non-trivial unitary 2-cocycle.\ethm
    %\begin{proof}
\blmma \label{lemma23}
For any choice of $(\tau,\lambda,k_0,\mu)$ and $\psi\equiv 1$
the corresponding fiber functor $\phi$ is monoidally isomorphic to the identity fiber functor.
\elmma
\begin{proof}
Let $v_{\chi}:\bbc_{\chi}\to\bbc_{\chi},v_{\pi_2}:\clh_{\pi_2}\to \clh_{\pi_2},v_{\pi'_2}:\clh_{\pi'_2}\to\clh_{\pi'_2},v_{\pi_4}:\clh_{\pi_4}\to\clh_{\pi_4}$ be the unitary linear maps given by
\begin{enumerate}
    \item [1)] $v_{\chi}(1_{\bbc_{\chi}}=)\tau(\chi)1_{\bbc_{\chi}},$
    \item[2)] $v_{\pi_2}=\lambda Id_{\clh_{\pi_2}},$
    \item[3)] $v_{\pi'_2}=\lambda k_0Id_{\clh_{\pi'_2}},$
    \item[4)] $v_{\pi_4}=\mu^{1/2}Id_{\pi_4}.$
\end{enumerate}
One can check that $(v_a\ot v_b)(\Upsilon^{a\ot b}_{c})v^{\ast}_c=\phi(\Upsilon^{a\ot b}_c)$ for any $a,b,c\in\{\chi_{ijk},\pi_2,\pi'_2,\pi_4).$ Hence $\phi$ corresponds to the identity tensor functor on $Corep(C(G)).$
\end{proof}

    \blmma \label{lemma24}
    When $\psi(\chi_{100})=1$ and $\psi_(\chi_{001})=-1,$ for any two choices of $(\tau_1,\lambda_1,k_01,\mu_1)$ and $(\tau_2,\lambda_2,k_{02},\mu_2)$ the corresponding fiber functor $\phi_1,\phi_2$ are monoidally isomorphic.
    \elmma
    \begin{proof}
    It is easy to observe that $\phi^{-1}_1\phi_2\cong Id.$ Hence $\phi_1$ is isomorphic to $\phi_2.$
    \end{proof}
    \blmma \label{25}
    Any two fiber functors as in lemma \ref{lemma23} and as in lemma \ref{lemma24}, are not monoidally isomorphic.
    \elmma
    \begin{proof}
        Without loss of generality, we can take first functor to be identity tensor functor and other to be $\phi$ corresponding $\lambda=1, k_0=1,\mu=1,\tau=1.$ Suppose $\phi\cong Id.$
         So there exists a unitary morphism $V_{\pi_2\ot\pi_4}\in Mor(\pi_2\ot\pi_4, \pi_2\ot\pi_4)$ such that $V_{\pi_2\ot \pi_4}(T_1)V^{\ast}_{\pi_4}=T_1, V_{\pi_2\ot \pi_4}(T_2)V^{\ast}_{\pi_4}=-T_2 $ where $T_{1},T_{2}\in Mor(\pi_4, \pi_2\ot\pi_4).$ As $\pi_4$ is irreducible, $V_{\pi_4}$ must be of the form $cI_{\pi_4}$ for some constant $c\in \bbc-\{0\}.$ Let $V_{\pi_2\ot\pi_4}=A.$ We already know $T_{1}(x_1)= e_1\otimes x_4$ and $T_2(x_1)=e_2\ot x_2$ . Hence $c^{-1}A(e_1\ot x_4)= e_1 \ot x_4$ and $c^{-1}A(e_2\ot x_2)=-e_2\ot x_2.$ This implies $c^{-1}A\neq Id$ and therefore $\phi$ is not isomorphic to the identity tensor functor.
    \end{proof}

As a corrollary, we get our final result
\bthm $H^{2}_{uinv}(C^{\ast}(G),S^{1})=Z_2.$
\ethm
\begin{proof}
The proof follows from lemmas \ref{lemma23}, \ref{lemma24}, \ref{25}. The only non-trivial class is given by any functor in lemma \ref{lemma24}.
\end{proof}

\brmrk 
Similarly we can prove that $H^{2}_{inv}(C^{\ast}(G), \bbc-\{0\})$ is $Z_2.$
\ermrk
   \section{Invariant 2-cocycles of dual of Kac-Paljutkin algebra}
Let us recall the Tambara-Yamagami tensor category \cite{TAMBARA1998692} .

  Tambara--Yamagami tensor categories \cite{TAMBARA1998692} is   equivalent to the
the category of representations of the Kac--Paljutkin Hopf algebra \cite{TAMBARA1998692}, which is arising from   the
Klein $4$-group $K_4 = \mathbb{Z}/2\mathbb{Z} \oplus \mathbb{Z}/2\mathbb{Z}.$
Elements  of $K_4 = \{ e,s,t,st \}$
satisfies the relations $s^2 = t^2= (st)^2 =e.$
% We focus on  the tensor category $\mathcal{C}(\chi, \tau)$ corresponding to 
$\chi = \chi_c$ is a nondegenerate symmetric bicharacter  of $K_4$
which is given by
\begin{align*}
\chi_c(a,a) = \chi_c(b,b) = -1,\ \chi_c(a,b) =1,
\end{align*}
and considering the parameter $\tau = \frac{1}{2}$ . Now, we define the category $\mathcal{C}(\chi, \tau)$ and
Its objects are finite direct sums of elements in $ S = K_4 \cup \{ \rho \}$.
Sets of morphisms between elements in $S$ are given by
\[
Mor (s,s') = 
\begin{cases}
	\mathbb{C} & s = s', \\
	0 & s \neq s',
\end{cases} 
\]
so $S$ is the set of irreducible classes of $\mathcal{C}(\chi, \tau)$.
Tensor products of elements in $S$ are given by  
\[
\ s \otimes \rho = \rho = \rho \otimes s, \ \rho \otimes \rho = \bigoplus_{s \in K_4} s, s \otimes t = st,\ (s,t \in K_4)
\]
and the unit object is $e$.
Associativities $\varphi$ are given by 
\begin{align*}
	\varphi_{s,t,u} &= \text{id}_{stu},&
	\varphi_{s,t,\rho} &= \varphi_{\rho,s,t} = \text{id}_{\rho},\\
	\varphi_{s,\rho,t} &=  \chi_c (s,t) \text{id}_{t}, &
	\varphi_{s,\rho,\rho} &= \varphi_{\rho,\rho,s} = \bigoplus_{k \in K_4} \text{id}_{k},\\
	\varphi_{\rho,s,\rho} &= \bigoplus_{k \in K_4} \chi_c (s,t) \text{id}_{k}, &
	\varphi_{\rho, \rho, \rho} &= \left( \frac{1}{2} \chi_c (k,l)^{-1} \text{id}_{\rho} \right)_{k,l} \colon \bigoplus_{k \in K_4} \rho \to \bigoplus_{l \in K_4} \rho,
\end{align*}
for $s,t,u \in K_4$. Now,  if we choose the natural fiber functor of this category then 
this category is identified with the corepresentation category of Kac--Paljutkin quantum group $\q_{kp}$, that is
\[
\mathcal{C} \left( \chi_c, \frac{1}{2} \right) \simeq Rep(\q_{kp})\simeq Corep(\hat{\q_{kp}} )
\]
as tensor categories.

Moreover, using the discussion and calculation in \cite{TAMBARA1998692}, we observe that there is a fiber functor $\phi_{0}$ from which $\q_{kp}$ is obtained by the Tannaka-Krein reconstruction. It can be seen from \cite{TAMBARA1998692} that $\phi_{0}(s)\cong\phi_{0}(t)\cong\phi_{o}(st)\cong \bbc$ and $\phi_{0}(\rho)=\bbc^{2}.$ Moreover we can choose the basis element $U_{s},U_{t},U_{st}$ and $V_{s},V_{t},V_{st}$:

 \begin{align*}
    U_{s}=\begin{pmatrix}
        0 & i\\
        1 & 0
    \end{pmatrix},
    U_{t}=\begin{pmatrix}
        0 & 1\\
        i & 0
    \end{pmatrix},
    U_{st}=\begin{pmatrix}
        -1& 0\\
        0 & 1
    \end{pmatrix},
\end{align*} where $U_{s}\in Mor(\rho, s \ot \rho), U_{t}\in Mor(\rho, t\ot \rho), U_{st}\in Mor (\rho, st\ot \rho)$ 
and
 \begin{align*}
    V_{s}=\begin{pmatrix}
        0 & 1\\
        i & 0
    \end{pmatrix},
    V_{t}=\begin{pmatrix}
        0 & i\\
        1 & 0
    \end{pmatrix},
    V_{st}=\begin{pmatrix}
        -1& 0\\
        0 & 1
    \end{pmatrix},
\end{align*} where $V_{s}\in Mor(\rho, \rho \ot s), V_{t}\in Mor(\rho, \rho\ot t)~~and~~ V_{st}\in Mor (\rho,  \rho\ot st).$
Now, let $\phi$ be a dimension preserving fiber functor on this category. 

Let $\phi(\Upsilon^{\rho\ot t}_{\rho}):=\widetilde{\Upsilon}^{{\rho\ot t}}_{\rho}$, $\phi(\Upsilon^{t\ot\rho}_{\rho}):=\widetilde{\Upsilon}^{t\ot\rho}_{\rho}$ and $\phi(\Upsilon^{\rho\ot \rho}_s):=\widetilde{\Upsilon}^{{\rho\ot \rho}}_s.$

Here we introduce some notations which are,
\begin{align*}
    &\widetilde{\Upsilon}^{{s\ot t}}_{st}=\theta(s,t)\Upsilon^{{s\ot t}}_{st},\\
   &\widetilde{\Upsilon}^{{\rho\ot t}}_{\rho}=d_t\Upsilon^{\rho\ot t}_{\rho}=d_tV_t,\\
   &\widetilde{\Upsilon}^{t\ot\rho}_{\rho}=c_t\Upsilon^{\rho\ot t}_{\rho}=c_tU_t,\\
   &\widetilde{\Upsilon}^{\rho\ot\rho}_{t}=k_t{\Upsilon}^{\rho\ot\rho}_{t}.
  \end{align*}

  Now if we choose the unitary linear maps $ v_1:\bbc_1\to \bbc_1, v_s:\bbc_s\to \bbc_s, v_t:\bbc_{t}\to \bbc_t,v_{st}:\bbc_{st}\to \bbc_{st} ~~ and ~~v_{\rho}:H_{\rho}\to H_{\rho}$ such that $v_{i}$ are identity maps from $\bbc_{i} ~to ~\bbc_{i}$ and $v_{\rho}=k_{1}^{-1/2}Id_{\clh_{\rho}}$ then it follows from the proof of Proposition (3.5)  \cite{etn} that $\phi$ is isomorphic to a fiber functor $\phi'$ where $\phi'(\Upsilon^{a\ot b}_c)=(v_a\otimes v_b)\phi(\Upsilon^{a\otimes b}_{c} )v^{\ast}_c$ for which $\phi'(\Upsilon^{\rho\ot\rho}_{e})=\Upsilon^{\rho\ot\rho}_{e}.$
Without loss of generality, Let us assume that $k_1=1.$
  \blmma $\theta$ is a 2-cycle on $K_4.$
  \elmma
\begin{proof} Proof of this lemma similar to the proof of lemma \eqref{le3.2}, hence omitted.
\end{proof}

Without loss of generality, we assume that $\theta$ is a normalized 2-cycle.

\blmma \label{le4.2}  $c_xc_y=\theta(x,y)c_{xy},$ where $x,y\in K_4.$
\elmma
\begin{proof} 
From the diagram 
\[\begin{tikzcd}
	{x\otimes y\otimes\rho} && {x\otimes \rho} \\
	{xy\otimes \rho} && \rho
	\arrow["{(Id_{\bbc_{x}}\otimes c^{\ast}_yU^{\ast}_{y})}", from=1-1, to=1-3]
	\arrow["{c^{\ast}_xU^{\ast}_x}"', from=1-3, to=2-3]
	\arrow["{\theta^{\ast}(x,y)\otimes Id_{H_{\rho}}}"', from=1-1, to=2-1]
	\arrow["{c^{\ast}_{xy}U^{\ast}_{xy}}", from=2-1, to=2-3],
\end{tikzcd}\] we can conclude that $c_xc_y=\theta(x,y)c_{xy}.$
\end{proof}
\blmma \label{le4.3} $d_x d_y=\theta(x,y)d_{xy},$ where $x,y\in K_4.$
\elmma
\begin{proof} Proof of this lemma similar to the previous lemma, hence omitted.
\end{proof}
\blmma$c_x=\tau(x)d_x,$ where $\tau$ is a character on $K_4.$
\elmma
\begin{proof}
From lemmas \eqref{le4.2} and \eqref{le4.3}, one can conclude that 
\begin{align}
    \theta (x,y)=\frac{c_xc_y}{c_{xy}}=\frac{d_{x}d_y}{d_{xy}}.
\end{align}
 Hence, $(c_xd^{-1}_x)(c_yd^{-1}_y)=c_{xy}d^{-1}_{xy}.$ This implies that $c_xd^{-1}_x=\tau(x)$ for a 2 cycle of $K_4.$
 \end{proof}

From the associativity relation, 
\begin{align*}
   \varphi_{s,\rho,t} &=  \chi(s,t) \text{id}_{t}, 
\end{align*} we observe that 
\begin{align}
    U_xV_y=\chi(x,y)V_yU_x.
\end{align}

Let $P_x$ be range of $\Upsilon ^{m\ot m}_x$ and assume that $\epsilon_x$ is the image of $P_x.$

We already know that $\varphi_{s,\rho,t}=\chi_{c}(s,t)id_t$. 

Now, one can easily observe from the associativity relations that 
\begin{align}
    U_xV_y=\chi(x,y)V_yU_x.
\end{align}
It is a straight forward computation to verify the following
\begin{align}
  &(U_s \ot Id_{\rho})\epsilon_s=i\epsilon_{1},\\
  &(U_t \ot Id_{\rho})\epsilon_t=i\epsilon_{1},\\
  &(U_{st}\ot Id_{\rho})\epsilon_{st}=(-1)\epsilon_{1}.
  \end{align}
 
\blmma \label{le4.5} The following identities hold: 
\begin{enumerate}
    \item[1)] $c_sk_t=\theta(s, s^{-1}t)k_{s^{-1}t},$
    \item[2)]$d_{s}k_t=\theta(s,s^{-1}t)k_{s^{-1}t}.$
\end{enumerate} 
\elmma
\begin{proof}
 $1)$ From this diagram
 \[\begin{tikzcd}
	{(s\ot\rho)\ot\rho} && {s\ot (\rho\ot\rho)} \\
	{\rho\ot\rho} && {s\ot s^{-1}t} \\
	& t
	\arrow["{\oplus_k Id_{k}}", from=1-1, to=1-3]
	\arrow["{c^{\ast}_sU^{\ast}_s\ot Id_{\rho}}"', from=1-1, to=2-1]
	\arrow["{Id_s\ot k^{\ast}_{s^{-1}t}P^{\ast}_{s^{-1}t}}", from=1-3, to=2-3]
	\arrow["{k^{\ast}_tP^{\ast}_t}"', from=2-1, to=3-2]
	\arrow["{\theta^{\ast}(s, s^{-1}t)}", from=2-3, to=3-2]
\end{tikzcd}\] 
one can easily conclude that  $c_sk_t=\theta(s, s^{-1}t)k_{s^{-1}t}.$

2) Similarly, we can prove that $d_{s}k_t=\theta(s,s^{-1}t)k_{s^{-1}t}.$

\end{proof}

\blmma   
$c_xk_x=1$ and  $c_x=d_x,$ where $x\in K_4.$
\elmma
\begin{proof}
 If we choose $s=t$ then it follows $c_sk_s=\theta(s,1)k_{1}=k_1=1$ from lemma   \eqref{le4.5}. Similarly, we can deduce that $d_sk_s=k_1$. Hence $c_s=d_s.$
\end{proof}
\bthm \label{thm4.7} $H^{2}_{uinv}(\hat{\q}_{kp}, S^{1})\cong1,$ $H^{2}_{inv}(\hat{\q}_{kp},\bbc-\{0\})\cong 1.$
\ethm
\begin{proof}
    We define unitary linear maps $v_1:\bbc_1\to \bbc_1, v_s:\bbc_s\to \bbc_s, v_t:\bbc_{t}\to \bbc_t,v_{st}:\bbc_{st}\to \bbc_{st} ~~ and ~~v_{\rho}:H_{\rho}\to H_{\rho}$, which are given by
\begin{align*}
    &v_1(1_{\bbc_1})=c_{1}1_{\bbc_1},\\
    &v_s(1_{\bbc_s})=c_{s}1_{\bbc_s},\\
    &v_t(1_{\bbc_t})=c_{t}1_{\bbc_t},\\
    &v_{st}(1_{\bbc_{st}})=c_{st}1_{\bbc_{st}},\\
    &v_{\rho}=Id_{H_{\rho}}.
\end{align*}

Now, one can check that $\phi(\Upsilon^{a\ot b}_c)= (v_a \ot v_b)(\Upsilon^{a\ot b}_c) v^{\ast}_c.$

Hence $H^{2}_{uinv}(\hat{\q}_{kp},S^{1})=1.$

Similarly, $H^{2}_{inv}(\hat{\q}_{kp}, \bbc-\{0\})\cong 1$ as $v^{\ast}_c=v^{-1}_c$ for $c\in\{1,s,t,\rho\}.$
\end{proof} 
\brmrk
In fact, it is well known that $\q_{kp}\cong \hat{\q}_{kp}$ as Hopf *-algebras, hence Theorem \eqref{thm4.7} is also valid if $\hat{\q}_{kp}$ is replaced by $\q_{kp}.$
\ermrk

\bibliographystyle{alpha}
\bibliography{ref}

\end{document}